\documentclass[final,onefignum,onetabnum]{siamart190516}

\usepackage{lipsum}
\usepackage{amsfonts}
\usepackage{graphicx}
\usepackage{epstopdf}
\usepackage{algorithmic}
\usepackage[utf8]{inputenc}
\usepackage{float}
\usepackage{amsmath}
\usepackage{xcolor}
\usepackage{appendix}
\usepackage[ruled,vlined,linesnumbered]{algorithm2e}

\ifpdf
  \DeclareGraphicsExtensions{.eps,.pdf,.png,.jpg}
\else
  \DeclareGraphicsExtensions{.eps}
\fi

\newcounter{todocounter}

\newsiamremark{remark}{Remark}
\newsiamremark{hypothesis}{Hypothesis}
\crefname{hypothesis}{Hypothesis}{Hypotheses}
\newsiamthm{claim}{Claim}

\headers{Differentiating SRB density on unstable manifolds}{Adam A. \'Sliwiak and Qiqi Wang}

\title{A trajectory-driven algorithm for differentiating SRB measures on unstable manifolds\thanks{Submitted to SIAM.
\funding{This work was funded by Air Force Office of Scientific Research Grant No. FA8650-19-C-2207 and U.S. Department of Energy Grant No. DE-FOA-0002068-0018.}}}

\author{Adam A. \'Sliwiak\thanks{Center for Computational Science and Engineering, Department of Aeronautics and Astronautics, Massachusetts Institute of Technology, Cambridge, Massachusetts, 02139, United States of America 
  (\email{asliwiak@mit.edu, qiqi@mit.edu}).}
\and Qiqi Wang\footnotemark[2]}

\usepackage{amsopn}

\begin{document}
\maketitle

% REQUIRED
\begin{abstract}
  SRB measures are limiting stationary distributions describing the statistical behavior of chaotic dynamical systems. Directional derivatives of SRB measure densities conditioned on unstable manifolds are critical in the sensitivity analysis of hyperbolic chaos. These derivatives, known as the {\it SRB density gradients}, are by-products of the regularization of Lebesgue integrals appearing in the original {\it linear response} expression. In this paper, we propose a novel trajectory-driven algorithm for computing the SRB density gradient defined for systems with high-dimensional unstable manifolds. We apply the concept of measure preservation together with the chain rule on smooth manifolds. Due to the recursive one-step nature of our derivations, the proposed procedure is memory-efficient and can be naturally integrated with existing Monte Carlo schemes widely used in computational chaotic dynamics. We numerically show the exponential convergence of our scheme, analyze the computational cost, and present its use in the context of Monte Carlo integration.
  
\end{abstract}

% REQUIRED

% REQUIRED

%This is for Elsevier article
%\begin{frontmatter}

%\journal{Journal Name}

%\title{A trajectory-based algorithm to compute SRB density gradient}
%\author{Adam A. \'Sliwiak\corref{cor}}
%\ead{asliwiak@mit.edu}
%\author{Qiqi Wang}
%\ead{qiqi@mit.edu}
%\address{Center for Computational Science and Engineering, Massachusetts Institute of Technology (MIT), 77 Massachusetts Avenue, Cambridge, MA, 02139, USA}
%\address{MIT Department of Aeronautics and Astronautics}
%\cortext[cor]{Corresponding author.}

%\begin{abstract}
%    This paper focuses on the computation of SRB density gradient for general higher-dimensional chaotic systems.
%\end{abstract}
%\begin{keyword}
%SRB measure, SRB density gradient, Ergodicity, Chaotic Dynamical Systems
%\end{keyword}
%\end{frontmatter}

\section{Introduction}\label{sec:introduction}
%What is SRB measure and why is it useful?
Due to their seemingly irregular and quasi-random behavior, a mathematical description of chaotic dynamical systems might be challenging. A major breakthrough in the analysis of chaos was the introduction of the SRB (Sinai-Ruelle-Bowen) measure $\mu$ \cite{ruelle-srb}. This scalar quantity, defined on a compact Riemannian manifold, contains a coherent statistical description of the dynamics. Intuitively, the SRB measure represents the likelihood of the trajectory passing through a non-zero-volume region of a strange attractor. Although the concept of SRB measures was originally applied to Axiom A systems, several rigorous studies extended this idea even beyond the universe of uniformly hyperbolic systems \cite{young-srb, climenhaga-srb, cruz-srb}. 

%What is SRB density gradient g and why is it useful?
Lebesgue integrals with respect to $\mu$, which represent expected values of certain smooth observables, are fundamental in the analysis of chaos. Under the assumption of {\em ergodicity}, they equal the time-average of an infinitely-long sequence generated along a trajectory. Integrals of this type can thus be approximated using a Monte Carlo method. If the integrand involves highly-oscillatory derivatives, then the Monte Carlo integration might be prohibitively expensive due to a large variance of the sample \cite{sliwiak-densitygrad}. In case of derivatives of functions evaluated at a future time (see examples of such integrands in \cite{chandramoorthy-s3-new, sliwiak-1d, gritsun-fdt, abramov-original}), the direct use of any integration scheme might be impossible due to the {\em butterfly effect}. Indeed, the application of the chain rule results in a product of the system's Jacobian matrices whose norms increase exponentially in time. A remedy for this computational difficulty is integration by parts, which moves the differentiation operator away from the problematic function to the SRB measure. This is in fact a consequence of the generalized fundamental theorem of calculus. In addition to the boundary term, we effectively obtain a new Lebesgue integral involving a product of the antiderivative of the original integrand and the {\em SRB density gradient} $g = \partial\log\rho =\partial \rho /\rho$, where $\rho$ denotes the density of $\mu$ (i.e., the Radon-Nikodym derivative \cite{nagy-density}).

%Appearance of g in literature. 
The SRB density gradient is critical in the sensitivity analysis of chaos. The major implication of Ruelle's {\em linear response} theory is a closed-form expression for the parametric derivative of long-time averages (a.k.a. the system's sensitivity) \cite{ruelle-original,ruelle-corrections}. The space-split sensitivity (S3) method \cite{chandramoorthy-s3, chandramoorthy-s3-new} reformulates Ruelle's formula to a computable form by splitting the perturbation vector and performing integration by parts on unstable manifolds. Using the S3 formula, one can construct an efficient and provably convergent Monte Carlo algorithm for sensitivities in uniformly hyperbolic systems. This algorithm requires computing the SRB density gradient defined as a directional derivative of $\rho$ conditioned on the unstable manifold. Indeed, the SRB measure is generally singular with respect to Lebesgue measure in the stable direction \cite{young-srb}. This approach of the regularization of Ruelle's integrals on unstable manifolds has also been applied in \cite{ni-fast} to derive a fast linear response algorithm for differentiating SRB states.   
Several algorithms for sensitivity analysis that stem from the Fluctuation-Dissipation Theorem (FDT) \cite{kubo-fdt} also require $g$ \cite{abramov-original,abramov-blended}. Motivated by empirical data of certain chaotic models, some FDT-based methods assume Gaussian distribution of measure \cite{gritsun-fdt}. Such an assumption reduces the FDT linear response operator to a simple time autocorrelation function, which dramatically facilitates the sensitivity computation for the cost of limited applicability. The density gradient can also be used as an reliable indicator of the differentiability of statistical quantities \cite{sliwiak-differentiability} in chaotic systems. In particular, the slope of the distribution tail of $g$ have been shown to be strictly associated with the existence of parametric derivatives of statistics. Therefore, we seek a numerical procedure for $g$ that does not make any assumptions about the statistical behavior of the system and is thus generalizable to any chaotic dynamical system that admits SRB measures. 

%Early attempts of computing g.
There already exist algorithms for the SRB density gradient derived for systems with one-dimensional unstable manifolds. In case of simple one-dimensional maps, one can derive an exponentially convergent recursion for $g$ using the measure preservation property \cite{sliwiak-1d}. The same formula can be inferred using the fact the SRB density is an eigenfunction of the Frobenius-Perron operator with eigenvalue 1 \cite{sliwiak-differentiability}. The authors of \cite{chandramoorthy-clv} propose an ergodic-averaging algorithm for self-derivatives (i.e., directional derivatives along one-dimensional expanding directions) of covariant Lyapunov vectors (CLVs) corresponding to the only positive Lyapunov exponent, which are tangent to unstable manifolds at any point on the attractor. Using the chain rule on smooth manifolds, one can show $g$ depends on the self-derivative of CLV at the previous time step, and this relation is governed by a second-order tangent equation \cite{chandramoorthy-clv, sliwiak-differentiability}. In a recent work, Ni \cite{ni-fast} proposes an algorithm for divergence on the unstable manifold using approximate {\it shadowing} coordinates instead of the full basis of the expanding subspace, as opposed to the S3 method.   

%What do we do in this paper?
In this paper, we systematically derive a trajectory-driven algorithm for the SRB density gradient by extending the measure preservation property to high-dimensional smooth manifolds. Using the density-based parameterization of unstable manifolds and the chain rule, it is possible to establish a recursive relation for the evolution of first- and second-order parametric derivatives of the coordinate chart. By definition, this chart is strictly associated with $g$ and can be interpreted as an SRB inverse cumulative distribution (quantile function). This type of parameterization, motivated by popular methods of statistical inference \cite{fiori-stats1}, has been thoroughly explained by the authors in \cite{sliwiak-densitygrad} in the context of simple Lebesgue measures. Through the relation of $g$, the coordinate map and its parametric derivatives, we show the density gradient can be computed by solving a collection of first- and second-order tangent equations. We also show that the recurring problem of the butterfly effect, which leads to exploding norms of tangent solutions, can be eliminated by iterative orthonormalization of the chart gradient. The major benefit of our derivation is that it is naturally translatable to a practicable algorithm that can be easily integrated with existing methods for sensitivity analysis of chaos.         

%Paper overview
This paper is structured as follows. In Section \ref{sec:preliminaries}, we introduce the SRB measure, its gradient, and highlight their importance in the field of chaotic dynamics. Subsequently, in Section \ref{sec:1D}, we apply the density-based parameterization for the description of unstable manifolds to derive recursive relations for the SRB density gradient. This derivation is followed by a numerical example involving a chaotic map with straight one-dimensional expanding subspaces. Section \ref{sec:general} generalizes all the concepts to high-dimensional chaotic maps with an arbitrary number of positive Lyapunov exponents (LEs). Based on the systematically derived iterative relations, a practicable algorithm for $g$ is thoroughly described. We analyze its cost, memory requirements and convergence. In the same section, we also demonstrate a numerical example of Monte Carlo integration, which requires the computation of $g$. Section \ref{sec:conclusions} summarizes this work.  

\section{SRB measure and its gradient: definitions and significance}\label{sec:preliminaries}
%Explain basic concepts: rigorously define SRB measure, ergodicity, etc.

Consider a diffeomorphic map $\varphi: M\to M$, $M\in\mathbb{R}^n$, $n\in\mathbb{Z}^+$ with an Axiom A attractor. Theorem 1 of \cite{young-srb} asserts that there exists an {\em invariant} and {\em physical} probability measure $\mu$ (and its density $\rho$), which satisfies:
\begin{enumerate}
    \item Invariance/conservation of measure condition:
    \begin{equation}
        \label{eqn:pre-invariance}
        \mu(A) = \mu(\varphi^{-1}(A))
    \end{equation}
    for any Borel subset $A\subset M$.
    \item Physicality condition: for any smooth $f:M\to\mathbb{R}$,
    \begin{equation}
        \label{eqn:pre-physicality}
        \int_{M}f(x)\,d\mu(x) = \int_{M} f(x)\,\rho(x)\,d\omega(x) = \lim_{N\to\infty}\frac{1}{N}\sum_{k=0}^{N-1}f\circ\varphi^{k}(x_0),
    \end{equation}
    where $d\omega$ denotes the Riemmanian volume element. The initial state $x_{0}$ is assumed to be $\mu$-typical and $\varphi^{k}(\cdot) = \varphi(\varphi^{k-1}(\cdot))$, $\varphi^{1} = \varphi$, $\varphi^{0} = \mathrm{Id}$.
    \item Absolute continuity: Conditional measure of $\mu$ denoted by $\tilde{\mu}_x$ and defined on the unstable manifold $U_x$ at point $x\in U_x$ is absolutely continuous (analogous property applies to the conditional density $\tilde{\rho}_x$).
    \item Singularity with respect to Lebesgue measure: $\mu$ is generally sharp in the stable direction (across unstable manifolds).
    \item Unit measure axiom (probability universe): $$\int_M d\mu(x) = \int_M \rho(x)\;d\omega(x) = \mu(M) = 1.$$
\end{enumerate}

\begin{figure}
    \centering
    \includegraphics[width = \textwidth]{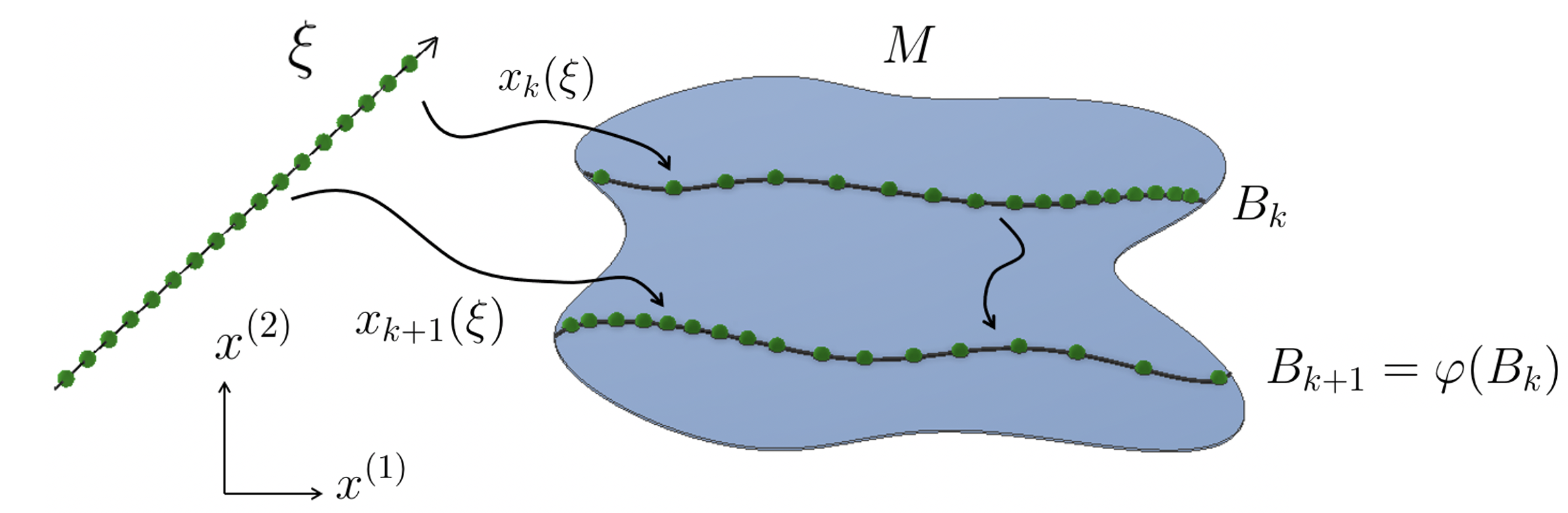}
    \caption{This figure graphically represents the measure preservation property. The localization of green bullets represents the SRB density on some 1D subspace of a 2D manifold $M$. All green bullets are equally weighted. In this sketch, we observe $\mu(B_k) = \mu(B_{k+1})$ and $B_{k+1}=\varphi(B_k)$, where $B_k\subset M$ and $B_{k+1}\subset M$ are parameterized by smooth charts, $x_k(\xi):[0,1]\to B_k$ and $x_{k+1}(\xi):[0,1]\to B_{k+1}$, respectively.}
    \label{fig:pre-properties}
\end{figure}

The measure $\mu$ and its density $\rho$ are respectively known as the SRB measure and SRB density distribution. We listed their properties most important in the context of this paper; however, the reader is referred to \cite{young-srb} for a detailed description of other significant features. One can think about Property (1) as the mass conservation law. For example, consider a Borel subset $B\subset{M}$ with a uniform measure that is mapped to $\varphi(B)\subset{M}$. If we divide $\varphi(B)$ into a finite number of subsets occupying the same volume, each of them generally has a different measure. In other words, each subset generally has its unique weight unless $\varphi$ represents a simple translation and/or rotation. Property (2) states that the SRB measure is physical, which means that by observing the system's evolution for an infinitely long period of time we can assign a weight (density) to each non-zero-volume region of the attractor. The expected value of any smooth function defined on $M$ can be computed as a simple volume integral over $M$ of that function multiplied by the density function. Figure \ref{fig:pre-properties} graphically explains Property (1), while the remaining four properties and their consequences are further explained and illustrated in the following sections. 

%Which dynamical systems have SRB measures? Need to add important comments -- need to expand on significance here
As mentioned above, SRB measures are guaranteed to exist in Axiom A (or, uniformly hyperbolic) systems. Different rigorous studies indicate that uniform hyperbolicity is in fact not required for the existence of $\mu$. For example, partially hyperbolic systems that have a mostly expanding \cite{alves-partial} or contracting \cite{burns-partial} central direction also admit SRB measures. In addition, many high-dimensional systems arising from discretization of real-world PDE models behave as uniformly hyperbolic systems, per the {\it hyperbolicity hypothesis} \cite{galavotti-hypothesis}.     

\begin{figure}
    \centering
    \includegraphics[width = \textwidth]{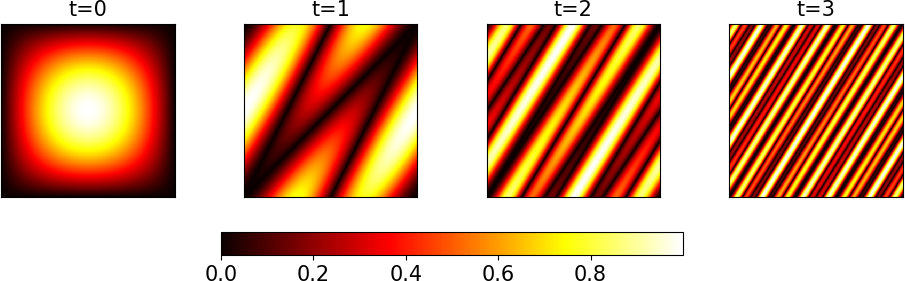}
    \caption{Evaluation of the composite function $f\circ\varphi^t(x)$ on the manifold $M=[0,1]^2$ at four consecutive steps $t$. In this case, the map $\varphi$ is the Arnold's cat map (Eq. \ref{eqn:pre-cat}), while $f(x^{(1)},x^{(2)}) = \sin(\pi x^{(1)})\,\sin(\pi x^{(2)})$. This particular $\varphi$ is a classical representative of an Anosov diffeomorphism.}
    \label{fig:pre-cat}
\end{figure}

%Significance of SRB measures
In many engineering applications, the expected value of some physically relevant quantity $f\in L^1(\rho)$, i.e., $\int_M f\,d\mu$, is usually of interest. The major challenge in the field of sensitivity analysis of chaos is to find a parametric derivative of the expected value, which is critical in grid adaptation \cite{larsson-grid}, optimization design \cite{jameson-conventional} and uncertainty quantification \cite{wang-phdthesis}. Ruelle rigorously derived a closed-form expression for that derivative \cite{ruelle-original,ruelle-corrections},
\begin{equation}
    \label{eqn:pre-sensitivity}
    \frac{d}{ds} \int_{M} f(x)\, d\mu(x) = \sum_{t=0}^{\infty}\int_{M} D(f\circ\varphi^t(x))\cdot \chi(x)\,d\mu(x),
\end{equation}
where $\chi$ denotes the derivative of $\varphi$ with respect to the map parameter $s$, while $D$ is a phase space differentiation operator. One could potentially apply a Monte Carlo algorithm to the integrals on the RHS of Eq. \ref{eqn:pre-sensitivity}. However, owing to the butterfly effect, the direct evaluation of the integrand for a higher $t$ is computationally infeasible. To illustrate this problem, let us consider the 2D Arnold's cat map $\varphi:[0,1]^2\to[0,1]^2$ defined as
\begin{equation}
    \label{eqn:pre-cat}
    x_{k+1} = Ax_{k}\,\text{mod}\,1,\;\;\;A=\begin{bmatrix}2 & 1 \\ 1 & 1\end{bmatrix},
\end{equation}
and some smooth function $f(x)$. In Figure \ref{fig:pre-cat}, we observe that even for a low $t$, $f\circ\varphi^{t}$ becomes highly-oscillatory, which implies that $\|D(f\circ\varphi^{t})\|$ grows very fast ($\|\cdot\|$ denotes the Euclidean norm in $\mathbb{R}^n$). Due to the presence of positive Lyapunov exponents in chaotic systems, the rate of growth is in fact exponential. It means that Ruelle's formula is impractical for a direct Monte Carlo computation.   

To circumvent this problem, one can apply integration by parts to move the differentiation operator away from the composite function. This idea was originally applied in the novel S3 method \cite{chandramoorthy-s3, chandramoorthy-s3-new}, and later also in Ni's approximate method \cite{ni-fast} for sensitivity analysis. In case of integrals with respect to a non-uniform measure, integration by parts requires differentiating the measure itself. However, according to Property (3) and Property (4), $D\rho$ generally does not exist. In this section, let us assume $\chi$ equals a unit vector $q$ that is tangent to the one-dimensional unstable manifold at every point on the manifold $M$\footnote{In a general case, $\chi\neq q$ and thus an extra step is required to regularize Ruelle's formula. This step involves a splitting of $\chi$ into two terms, such that one term belongs to unstable manifolds everywhere on $M$. The reader is referred to \cite{chandramoorthy-s3, chandramoorthy-s3-new} for a detailed description of this process.}. Thus, every integral from the RHS of Eq. \ref{eqn:pre-sensitivity} can be regularized through partial integration as follows,
\begin{align}
    \label{eqn:pre-integral-lhs}I &= \int_M Df_t(x)\cdot q(x)\,d\mu(x) \\
    \label{eqn:pre-integral-desint}  &= \int_{M/U}\int_{U_x} Df_t(s)\cdot q(s)\,d\tilde{\mu}_x(s)\,d\hat{\mu}(x) \\
    \label{eqn:pre-integral-dens}  &= \int_{M/U}\int_{U_x} \partial_q f_t(s)\,\tilde{\rho}_x(s)\,ds\,d\hat{\mu}(x) \\
    \label{eqn:pre-integral-param}  &= \int_{M/U}\int_{0}^1 \partial_{\xi} f_t(s(\xi))\,\tilde{\rho}_x(s(\xi))\,d\xi\,d\hat{\mu}(x) \\
    \label{eqn:pre-integral-parts}  &= -\int_{M/U}\int_{U_x} f_t(s)\,\partial_q\tilde{\rho}_x(s)\,ds\,d\hat{\mu}(x)+(\text{boundary term}) \\
    \label{eqn:pre-integral-reshuffle}  &= -\int_{M/U}\int_{U_x} f_t(s)\,\frac{\partial_q\tilde{\rho}_x}{\tilde{\rho}_x}(s)\,d\tilde{\mu}_x(s)\,d\hat{\mu}(x)+(\text{boundary term}) \\
    \label{eqn:pre-integral-rhs}  &= -\int_M f_t(x)\,g(x)\,d\mu(x) +(\text{boundary term}),
\end{align}
where $f_t(x):=f\circ(\varphi^t(x))$. To derive the final form of $I$, we perform the following steps. First (Step \ref{eqn:pre-integral-desint}), we disintegrate $\mu$ on a measurable partition $U$ determined by the geometry of unstable manifolds. The quotient measure $\hat{\mu}$ is defined such that for all Borel sets $B\subset M$, $$\mu(B) = \int_{M/U}\tilde{\mu}_x(B\cap U_x)\,d\hat{\mu}(x),$$
where $\tilde{\mu}_x$ is a conditional SRB measure with density $\tilde{\rho}_x$. Subsequently, in Step \ref{eqn:pre-integral-dens}, we use the measure-density relation, $d\tilde{\mu}_x = \tilde{\rho}_x\,ds$, where $s$ denotes the path length as we move along $U_x$. In Step \ref{eqn:pre-integral-param}, we parameterize $U_x$, which gives rise to $ds = \|x'(\xi)\|\,d\xi$. Note the multiplicative factor is absorbed by the parametric derivative of $f$, because $\partial_{\xi}f = \|x'(\xi)\|\,\partial_s f$. Integration by parts is applied in Step \ref{eqn:pre-integral-parts}, where the differentiation operator is moved from $f$ to $\tilde{\rho}$. In Steps \ref{eqn:pre-integral-reshuffle}-\ref{eqn:pre-integral-rhs}, we reshuffle terms and use the above identities again to simplify the final expression. Integration by parts also gives rise to a boundary term, which involves two integrals with respect to the quotient measure of $f\,\tilde{\rho}_x$ evaluated at $\xi = 0$ and $\xi = 1$, respectively. From now on, we shall drop the subscript notation for conditional distributions; the tilde $\tilde{(\cdot)}$ notation shall imply the given distribution is restricted to a local unstable manifold. Note the boundary term,
\begin{equation}
    \label{eqn:pre-bt}
    (\text{boundary term}) = \int_{M/U}\, \left[\tilde{\rho}(\xi)\,f_k(\xi)\right]_{\xi=0}^{\xi=1}\,d\hat{\mu}(x),
\end{equation}
can be expressed in terms of a regular volume integral over $M$ of the divergence on unstable manifolds,
which vanishes according to Theorem 3.1(b) of \cite{ruelle-original}. This is indeed a direct consequence of the fact the boundary terms across two neighboring rectangles of the Markov partition of $M$ cancel out. To visualize this property, let us consider the Arnold's cat map (Eq. \ref{eqn:pre-cat}), for example. Despite its ``artificial" discontinuities due to the modulo operator, this nonlinear transformation in fact maps a smooth torus to itself. One could arbitrarily change the boundaries of the square $M$ in both phase space directions without modifying the map itself, and still describe the same torus.    

Since $\int_M Df_k\cdot q\,d\mu = -\int_M f_k\,g\,d\mu$, we can thus alternatively apply Monte Carlo to the RHS that involves the SRB density gradient $g$ \cite{chandramoorthy-s3-new, sliwiak-differentiability, sliwiak-1d},
\begin{equation}
    \label{eqn:pre-g-def}
    g(x) = \frac{\partial_q\rho(x)}{\rho(x)}  = \frac{\partial_q\tilde{\rho}(x)}{\tilde{\rho}(x)} = \partial_q\log\tilde{\rho}(x).
\end{equation}
Note that the integrand appearing in the regularized version of $I$ does not grow exponentially with $t$ if $f$ is bounded, which makes the sensitivity formula computable (immune to the butterfly effect). The integration by parts, as presented above, is generally useful if the integrand involves highly-oscillatory functions. The Monte Carlo integral example presented in \cite{sliwiak-densitygrad} shows that the partial integration may reduce the number of samples a few orders of magnitude to achieve the desired approximation error. Therefore, the computation of $g$ might be beneficial not only in the context of the Ruelle/S3/FDT-based method for sensitivity approximation, but also in a general setting when the expected value of an ill-behaved quantity of interest in a chaotic system is needed. The following two sections focus on the computation of $g$ for systems with an arbitrary number of positive LEs. The primary goal is to derive a recursive procedure compatible with Monte Carlo algorithms, which are widely used in the field. 

\section{Computing SRB density gradient for systems with one-dimensional unstable manifolds}\label{sec:1D}

In this section, we consider a generic $n$-dimensional, $n\in\mathbb{Z}^{+}$, uniformly hyperbolic dynamical system with one-dimensional unstable manifold governed by the $C^2$ diffeomorphic map $\varphi: M\to M$. $M$ is thus a Riemannian manifold immersed in $\mathbb{R}^n$. There exists a measurable partition $U$ of $M$ such that each member of that partition, $U_x$, coincides with the unstable manifold that contains $x\in{M}$. In this particular case, each $U_x\subseteq M$ is geometrically represented by a curve embedded in $\mathbb{R}^{n}$. We strive to compute the directional derivative of the logarithmic SRB density $g$ defined by Eq. \ref{eqn:pre-g-def}.
%Add a one-sentence summary of the next two sections

\subsection{Derivation of the iterative formula}\label{sec:1D-derivation}

The following notation is used throughout this section. Let $x_k(\xi):[0,1]\to U_{k}\subset M$ denote a $C^2$ chart (diffeomorphic map) that describes the unstable manifold $U_{k}$, $k\in\mathbb{Z}$. For any $k$, two charts $x_k (\xi)$ and $x_{k+1}(\xi)$, defined respectively on $U_{k}$ and $U_{k+1}$, are related as follows,
\begin{equation}
    \label{eqn:1D-map}
    x_{k+1}(\xi) = \varphi(x_k(\xi))
\end{equation}
for all $\xi\in[0,1]$ (see Figure \ref{fig:pre-properties} for an illustration of an $n=2$ case). We use $D\varphi$ and $D^2\varphi$ to respectively denote the Jacobian ($n\times n$ matrix) and Hessian ($n\times n\times n$ third-order tensor) of $\varphi$.  Since $\varphi$ is invertible, Eq. \ref{eqn:1D-map} can be viewed as a mathematical description of the evolution of SRB measure. For any observable $f$ defined on $M$, evaluated along a certain trajectory, we use the following short-hand notation, $f\circ x_{k}(\xi):=f_k$. Derivatives of the chart with respect to the parameter $\xi$ are denoted using the prime ($'$) symbol. A reference to the $i$-th component of a vector/matrix/tensor is indicated inside round brackets located in the superscript; for example, $q^{(i)}$ denotes the $i$-th component of $q$. Finally, we use $\partial_i$ to denote differentiation with respect to the $i$-th coordinate of phase space.

Let us parameterize $U_{k}$ such that
\begin{equation}
\label{eqn:1D-int}
    \xi = \int_{\mathcal{C}_{{k}}(\xi)} \tilde{\rho}(x_{k}(\xi))\,ds,
\end{equation}
where $\mathcal{C}_{{k}}(\xi)$ represents the segment of $U_{k}$ between $x_k(0)$ and $x_k(\xi)$, which implies that $\mathcal{C}_{k}(1) \equiv U_{k}$. Consequently, 
$\tilde{\rho}_k$ is the conditional SRB density restricted to $U_{k}$ satisfying $\tilde{\rho}_k = \rho_k / \int_{U_{k}} \rho_k \,ds$. We call it a {\em measure-based} parameterization, as the value of the parameter $\xi$ coincides with the value of SRB measure at $x_k(\xi)\in U_{k}$. The variable transformation between $\xi$ and the arc length $s$ implies that
\begin{equation}
    \label{eqn:1D-trans}
    \tilde{\rho}(x_k(\xi))\,\|x'_k(\xi)\| = 1.
\end{equation}
Note Eq. \ref{eqn:1D-trans} is in fact a formula for the density change from a uniform to nonuniform distribution due to the nonlinear variable transformation $x_k(\xi)$. Since $\xi\in[0,1]$, $\tilde{\rho}_k$, $\xi$, $x_k$ can be respectively viewed as a probability density function (PDF), cumulative distribution (CDF), and inverse cumulative distribution (quantile function).
Using the measure-based parametrization as described above, the SRB density gradient, defined in Eq. \ref{eqn:pre-g-def}, can be expressed in terms of parametric derivatives of the chart, i.e.,
\begin{equation}
\label{eqn:1D-g}
    g(x_k(\xi))=\frac{\partial_{q}\tilde{\rho}_k}{\tilde{\rho}_k}(x_k(\xi)) = -\frac{x_k'(\xi)\cdot x_k''(\xi)}{\|x_k'(\xi)\|^3},
\end{equation}
for any $\xi\in[0,1]$. Here, the derivative $\partial_q$ is computed in the direction of increasing value of $\xi$. The reader is referred to the authors' previous work in \cite{sliwiak-densitygrad}, where Eq. \ref{eqn:1D-g} is derived by differentiating Eq. \ref{eqn:1D-trans}, and comprehensively described using various numerical examples.

We notice $x_k'(\xi) = \|x_k'(\xi)\|\;q(x_k(\xi))$, and rewrite Eq. \ref{eqn:1D-g} to
\begin{equation}
\label{eqn:1D-gk}
    g(x_k(\xi)) = - q(x_k(\xi))\cdot\frac{x_k''(\xi)}{\|x_k'(\xi)\|^2} := - q(x_k(\xi))\cdot a(x_k(\xi)) = -q_k\cdot a_k. 
\end{equation}
Eq. \ref{eqn:1D-gk} indicates that the magnitude of the SRB density gradient equals the length of the projection of the (re-scaled) curve acceleration vector on the line tangent to the curve.
We now use Eq. \ref{eqn:1D-map}, differentiate it twice with respect to $\xi$, and apply the chain rule to obtain the following expression,
\begin{equation}
\label{eqn:1D-chain}
x''_{k+1}(\xi) = D^2\varphi(x_k(\xi))(x_k'(\xi), x_k'(\xi)) + D\varphi(x_k(\xi))\;x_k''(\xi),
\end{equation}
which means that
\begin{equation}
\label{eqn:1D-a}
\begin{split}
    a(x_{k+1}(\xi)) = & \frac{x''_{k+1}(\xi)}{\|x'_{k+1}(\xi)\|^2} = \\ &
    \frac{\|x'_{k}(\xi)\|^2 D^2\varphi(x_k(\xi))\left( q(x_k(\xi)),q(x_k(\xi))\right)}{\|x'_{k+1}(\xi)\|^2} + \frac{D\varphi(x_k(\xi)) x''_k(\xi)}{\|x'_{k+1}(\xi)\|^2}.     \end{split}
\end{equation}
The bilinear form that appears in the first term on the RHS of Eq. \ref{eqn:1D-a} can be expressed using Einstein's summation convention, i.e., $[D^2\varphi(q,q)]^{(ijk)} = \partial_i\partial_j\varphi^{(k)}\,q^{(i)}\,q^{(j)}$.

Given $\|x'_{k+1}(\xi)\| = \alpha(x_{k}(\xi)) \|x'_{k}(\xi)\|$, where $\alpha(x_{k}(\xi)) = \|D\varphi(x_{k})\;q(x_k(\xi))\|$, we have
\begin{equation}
\label{eqn:1D-ak}
    a_{k+1} = \frac{(D^2\varphi)_k (q_k, q_k) + (D\varphi)_k\;a_k}{\alpha_k^2}.
\end{equation}
From the parametric derivative of Eq. \ref{eqn:1D-map} and the definition of $\alpha(x_{k}(\xi))$, the recursion
\begin{equation}
\label{eqn:1D-qk}
    q_{k+1} = \frac{D\varphi_k\;q_k}{\alpha_k}
\end{equation}
automatically follows. We emphasize the fact the above procedure for $g$ (involving Eq. \ref{eqn:1D-gk}, \ref{eqn:1D-ak}, \ref{eqn:1D-qk}) is completely analogous to the algorithm proposed in Section 4.2 of \cite{sliwiak-densitygrad}, which was meant for simple Lebesgue measures evolving due to a generic non-chaotic diffeomorphisms. Here, however, we consider the evolution of the SRB measure in a chaotic system. Due to the butterfly effect, the tangent solution exponentially increases in norm. Therefore, we need the normalizing factor $\alpha$ in the iterative formula for $a$ and $q$ along the trajectory. Since $\varphi$ is uniformly hyperbolic, the solution to the tangent equation in Eq. \ref{eqn:1D-qk} converges exponentially in $k$ to the backward Lyapunov vector that is tangent to the unstable manifold regardless of the choice of the initial condition $q_0$. Under the same assumption, the recursion in Eq. \ref{eqn:1D-ak} for the acceleration vector $a$ also converges uniformly in $k$ at an exponential rate to the true solution for any initial condition $a_0$ bounded in norm. The reader is referred to Lemma 7.7 in \cite{chandramoorthy-s3-new} for the proof of the preceding statement.

To summarize, using the measure-based manifold parameterization, we derived a simple recursive procedure for the SRB density gradient that exponentially converges in case of uniformly hyperbolic systems and does not depend on initial conditions. As for now, we restrict ourselves to systems with one-dimensional unstable manifolds. Our main intention here is to introduce basic concepts before we move to general cases in Section \ref{sec:general}.   

\subsection{Numerical example: computing SRB density gradient on straight unstable manifolds}\label{sec:1D-example1}

As a pedagogical example, let us consider a family of $n$-dimensional maps, $n\in\mathbb{Z}^+$, whose unstable manifolds are straight and, without loss of generality, aligned with the first coordinate of the phase space. Certainly, this family includes, but is not limited to, all one-dimensional chaotic maps. In this particular case, $q^{(i)} = \delta^{(i1)}$, where $\delta$ denotes the Kronecker delta. Consequently, the parametric derivative of the chart $x_k(\xi)$, for any $k$, has all zero entries except the first one and, therefore, $\alpha(x_k(\xi)) = |\partial_1\varphi^{(1)}(x_k(\xi))|$. Thus, our recursive algorithm for $g$, which involves Eq. \ref{eqn:1D-gk}, \ref{eqn:1D-ak}, and \ref{eqn:1D-qk}, reduces to a single scalar iterative formula,
\begin{equation}
    \label{eqn:straight-simple}
    g(x_{k+1}(\xi)) = \frac{g(x_{k}(\xi))}{\partial_1\varphi^{(1)}(x_k(\xi))} - \frac{\partial^2_1\varphi^{(1)}(x_k(\xi))}{(\partial_1\varphi^{(1)}(x_k(\xi)))^2}
\end{equation}
for all $\xi\in[0,1]$. We were allowed to drop the absolute values, because $x'_k(\xi) > 0$, which is a consequence of our choice of the manifold parameterization. In this simple event of a straight unstable manifold, only two scalars are required to advance the iteration, i.e., first- and second-order derivative (in phase space) of the first component of $\varphi$, since the map is expanding only in one direction. This result is fully consistent with early non-systematic attempts to construct such a procedure for $g$ in \cite{sliwiak-1d, sliwiak-differentiability}. The previous studies used the measure preservation property to derive an analogous version of Eq. \ref{eqn:straight-simple}.

To verify the correctness of our procedure, we consider the 2D perturbed Baker's map $\varphi: M\to M$, with $M=[0,2\pi]^2$, defined as follows \cite{chandramoorthy-s3-new},
\begin{equation}
\label{eqn:straight-bakers}
\begin{split}
    x_{k+1} = \varphi(x_k) = &\left(\begin{bmatrix}2x_k^{(1)}    \\ x_k^{(2)}/2 + \pi \lfloor x_k/\pi\rfloor\end{bmatrix} \right.+ \\ & \left.\begin{bmatrix}s_1/2\,\sin(x_k^{(1)}/2)+s_{2}/2\,\sin(2x_k^{(1)})\,\sin(x_k^{(2)})\\s_3\,\sin(x_k^{(2)})+s_4/2\,\sin(2x_k^{(1)})\,\sin(x_k^{(2)})\end{bmatrix}\right)\,\mathrm{mod}\,2\pi,
\end{split}
\end{equation}
where $s_{1},s_{2}, s_{3},s_{4}$ are real-valued map parameters. If all of them are zero, we obtain the classical Baker's map (first term of the RHS of Eq. \ref{eqn:straight-bakers}), which is named after the kneading operation that bakers apply to a two-dimensional square dough. In particular, the dough is first stretched horizontally (in the unstable direction) by a constant factor, then compressed vertically (in the stable direction) by the same factor, and so forth. The square-shaped domain is stretched to a $2\times 1$ rectangle, cut into two squares, which are subsequently stacked horizontally. The Baker's map is an invertible chaotic map with one positive and one negative Lyapunov exponent. 

By introducing an extra term proportional to the four parameters, we perturb the kneading operation in the direction not necessarily aligned with the phase space directions. Indeed, by manipulating these parameters' values, we can control the shape of the unstable manifold, which gives us an excellent study case in the context of the SRB gradient computation. Notice, for example, if $s_{4} = 0$ and $s_{3}$ is sufficiently small, the iteration in Eq. \ref{eqn:1D-qk} produces $q_k$ whose second coordinate, $q_k^{(2)}$, converges exponentially to zero with $k$. In this case, therefore, unstable manifolds are straight and aligned with the $x^{(1)}$-axis. We use this observation to design our first numerical test.

\begin{figure}
    \centering
    \includegraphics[width=0.47\linewidth]{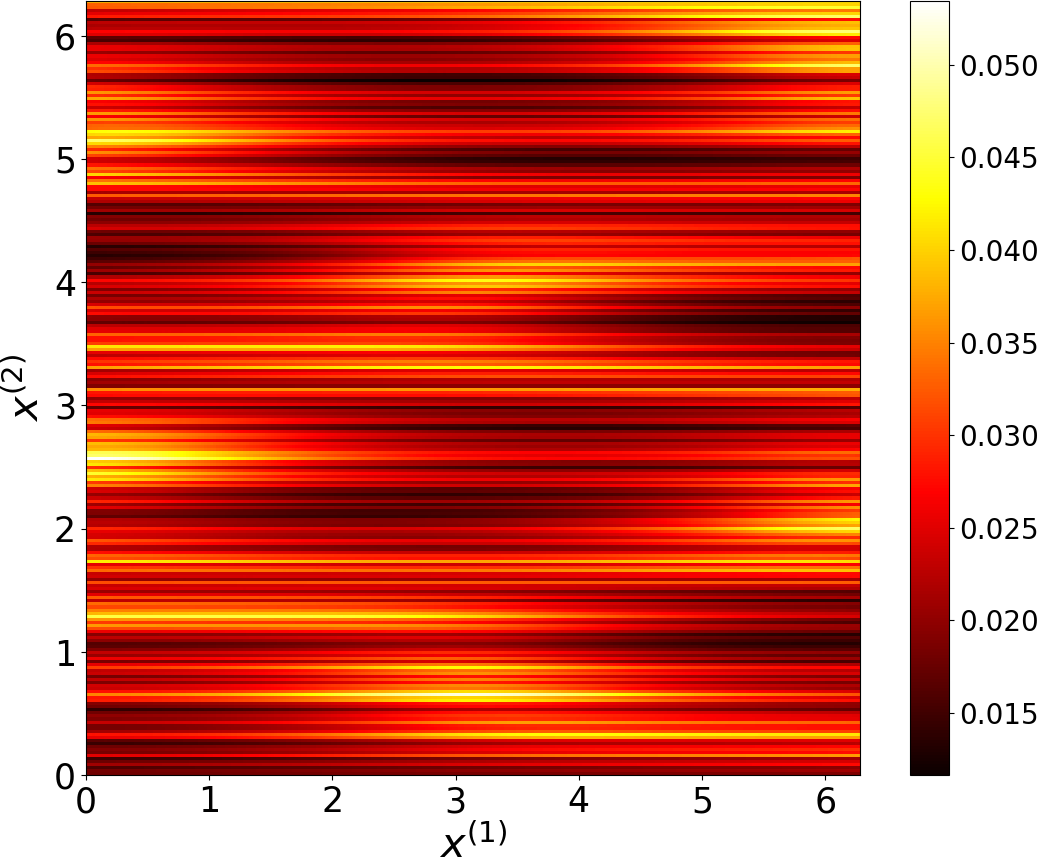}
    \includegraphics[width=0.47\linewidth]{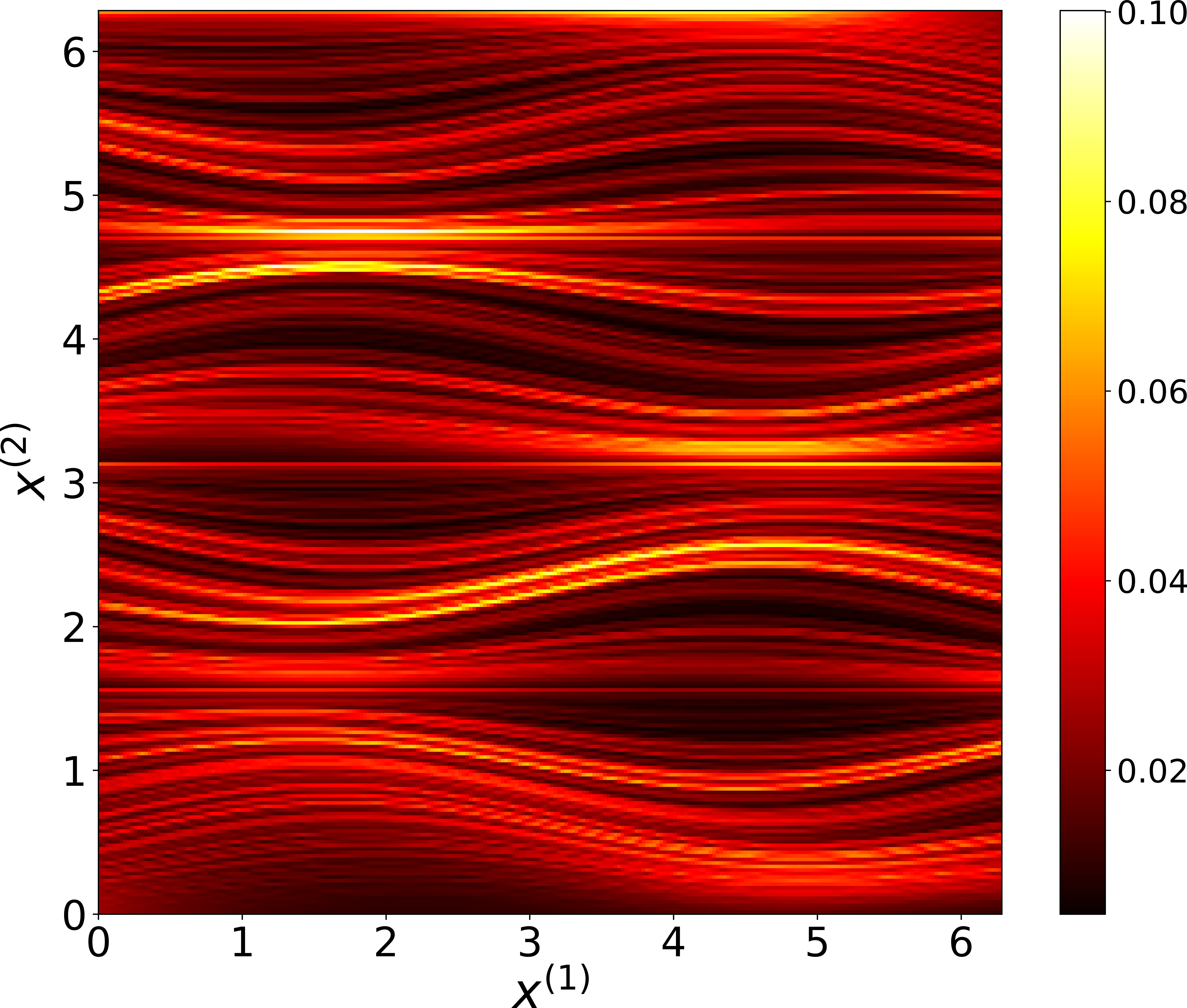}
    \caption{SRB distribution of the Baker's map with $s_{1}=s_{3}=s_{4}=0$, $s_{2} = 0.4$ (left plot) and $s_{1}=s_{2}=s_{3}=0$, $s_{4} = 0.4$ (right plot). We divided $M$ into $256^2$ rectangular bins of equal width and counted the number of times the trajectory passed through each of these bins. In this experiment, we generated $8000$ trajectories of length $209,715,200$, which gives us a total of approximately $1.68\cdot 10^{12}$ samples.}
    \label{fig:bakers_straight_pic}
\end{figure}

\begin{figure}
    \centering
    \includegraphics[width=0.46\linewidth]{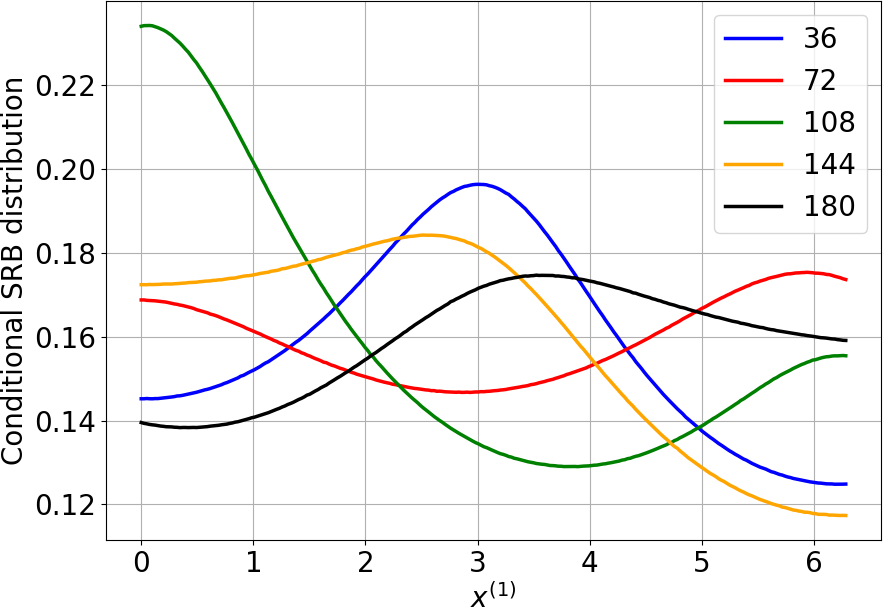}
    \includegraphics[width=0.47\linewidth]{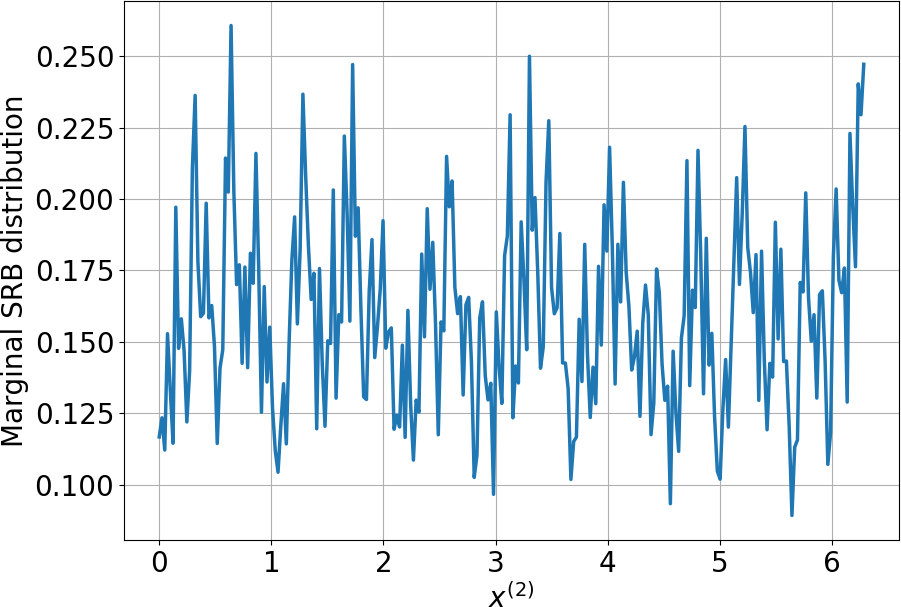}
    \includegraphics[width=0.6\linewidth]{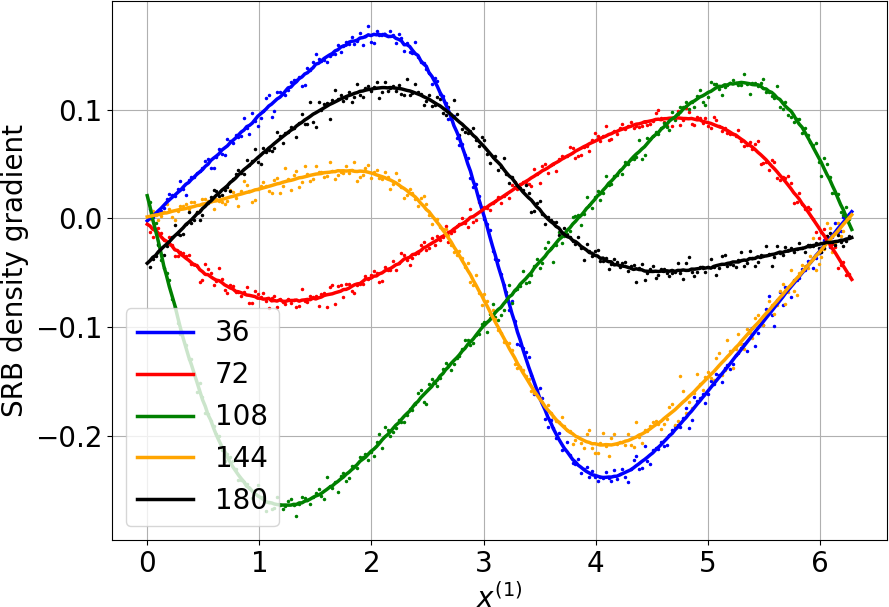}
    \caption{Upper left plot: conditional SRB distributions (SRB densities) corresponding to five different unstable manifolds. The numbers $36,72,108,144,180$ appearing in the legend represent the index of the horizontal bin row. For example, the red line corresponds to the SRB density defined on the unstable manifold at $x^{(2)}\approx 72/256 \cdot 2\pi \approx 1.76$. Upper right plot: marginal SRB distribution obtained through integrating the first coordinate out. Lower plot: SRB density gradient $g$ corresponding to SRB densities plotted in Figure \ref{fig:bakers_straight_conditional_marginal}. The $g$ function was computed using two distinct approaches: through the simplified trajectory-based recursion (Eq. \ref{eqn:straight-simple})(solid lines), and the central finite-difference method (dots). The oscillation of the finite-difference approximation is a manifestation of the statistical noise.}
    \label{fig:bakers_straight_conditional_marginal}
\end{figure}

In the first experiment, we consider the Baker's map defined by Eq. \ref{eqn:straight-bakers} with $s_{1}=s_{3}=s_{4}=0$ and $s_{2} = 0.4$. The left-hand side plot in Figure \ref{fig:bakers_straight_pic} illustrates the normalized SRB distribution corresponding to this parameter choice, which represents the probability of the trajectory passing through each square bin everywhere on $M$ (see the caption of Figure \ref{fig:bakers_straight_pic} for more details; for completeness, we also included a case with $s_{4} \neq 0$).
We observe a smooth behavior of the SRB distribution with respect to $x^{(1)}$ at any vertical level $x^{(2)}$. However, as we travel vertically, in the stable direction, the SRB distribution varies sharply. These radically different behaviors are typical symptoms of Property (3) and Property (4) of $\mu$ described in Section \ref{sec:preliminaries}, and they can also be observed in Figure \ref{fig:bakers_straight_conditional_marginal}, where the conditional and marginal SRB distributions are plotted, using data from Figure \ref{fig:bakers_straight_pic}. 

In Figure \ref{fig:bakers_straight_conditional_marginal}, we also plot the SRB density gradients defined on five different unstable manifolds. To compute $g$, the simplified recursion from Eq. \ref{eqn:straight-simple} was directly applied. To validate our computation, we approximated $g$ by applying the central finite-difference method to SRB densities plotted above. We observe a good agreement between the results obtained with these two different approaches, which confirms the correctness of our algorithm.

To conclude, in case of straight unstable manifolds, the SRB density gradient can be computed using the simplified recursive relation along trajectory (Eq. \ref{eqn:straight-simple}), which we verify through finite-differencing. This iteration is computationally cheap, as it involves solving a scalar tangent equation featuring both the first and second derivative of the first component of $\varphi$. In Appendix \ref{app:noninjective}, we show Eq. \ref{eqn:straight-simple} can also be applied to popular one-dimensional maps that are non-injective. We argue that certain non-measure-preserving  transformations have their higher-dimensional analogs similar to the classical Baker's map. Appendix \ref{app:hyperbolicity} presents a numerical study confirming the hyperbolicity of the Baker's map.

\section{Computing SRB density gradient for systems with general unstable manifolds}\label{sec:general}
We shall generalize the concepts introduced in Section \ref{sec:1D} to systems with $m$-dimensional unstable manifolds, $m\in\mathbb{Z}^+$. In other words, we consider general $n$-dimensional chaotic systems that have $m$ positive LEs, $1\leq m\leq n$. In this setting, the chart $x_k(\xi)$, $k\in\mathbb{Z}^+$, is a diffeomorphism that maps an $m$-dimensional hypercube, $[0,1]^m$, to the local unstable manifold $U_k\subset M$. For example, if $m = 2$ and $n = 3$, then the system has two positive LEs and its unstable manifolds are planes immersed in $\mathbb{R}^3\supset M$.

\subsection{Derivation of the iterative formula}\label{sec:general-derivation}

As introduced above, let us consider an $m$-dimensional smooth unstable manifold $U_k$ described by the chart $x_k(\xi):[0,1]^m\to U_k\subset M$. The vectors $x_k = [x_k^{(1)},...,x_k^{(n)}]^T$ and $\xi = [\xi^{(1)},...,\xi^{(m)}]^T$ have $n$ and $m$ components, respectively, and $0\leq \xi^{(i)}\leq 1$, $i=1,...,m$. We use $\nabla_\xi x_k(\xi)$ to denote the parametric gradient tensor of the chart. The $i$-th column of $\nabla_\xi x_k(\xi)$ contains the derivative of $x_k(\xi)$ with respect to $\xi^{(i)}$, i.e., $\partial_{\xi^{(i)}}x_k(\xi)$. For any Borel subset $V\subset [0,1]^m$ such that $x_k(V) = B_k \subset U_k$, the SRB measure-density relation can be expressed as follows,
\begin{equation}
    \mu(V) = \int_{B_k}\tilde{\rho}_k(x)\,d\omega(x),
\end{equation}
where $d\omega(x)$ denotes the natural volume element defined everywhere on $U_k$. Analogously to the 1D case described in Section \ref{sec:1D}, $\tilde{\rho}_k$ represents the conditional SRB density defined on $U_k$.
If we QR-factorize the parametric gradient of $x_k(\xi)$, 
\begin{equation}
    \label{eqn:general-qr}
    \nabla_{\xi} x_k(\xi) = Q(x_k(\xi))\;R(x_k(\xi))
\end{equation}
at any $\xi\in [0,1]^m$, the density conservation property can be expressed as
\begin{equation}
    \label{eqn:general-conserv}
    \tilde{\rho}(x_k(\xi))\;|\det R(x_k(\xi))| = 1,
\end{equation}
which is a generalization of Eq. \ref{eqn:1D-trans}. By differentiating Eq. \ref{eqn:general-conserv} with respect to $\xi$ and applying a non-trivial chain rule, we obtain
\begin{equation}
    \label{eqn:general-g2}
\begin{split}
    g^{(i)}(x_k(\xi)) =& \partial_{Q^{(:i)}(x_k(\xi))}\log\tilde{\rho}(x_k(\xi)) =\\ &-\frac{\mathrm{tr}\left(Q^T(x_k(\xi)) \;\partial_{\xi^{(i)}}\nabla_{\xi}x_k(\xi)\;R^{-1}(x_k(\xi))\right)}{\|\partial_{\xi^{(i)}}x_k(\xi)\|},
\end{split}
\end{equation}
or, equivalently,
\begin{equation}
\label{eqn:general-g2-einstein}
g^{(i)}(x_k(\xi)) = -\frac{Q^{(:j)}(x_k(\xi))\cdot\partial_{
\xi^{(i)}} \partial_{\xi^{(k)}} x_k(\xi)\;(R^{-1})^{(kj)}(x(\xi))}{\|\partial_{\xi^{(i)}}x_k(\xi)\|}
\end{equation}
for all $\xi\in [0,1]^m$, where the repeated indices imply summation (Einstein's convention), while the superscript $:i$ denotes the $i$-th column of a matrix. This expression was obtained by employing the orthonogonality of $Q$ and upper-triangular structure of $R$. It is computationally convenient as it does not involve parametric derivatives of the determinant of $R$. The reader is referred to \cite{sliwiak-densitygrad} for a step-by-step derivation of Eq. \ref{eqn:general-conserv}-\ref{eqn:general-g2-einstein}.

The purpose of this section is to derive an iterative (trajectory-driven) procedure for $g$. Analogously to the derivation in Section \ref{sec:1D}, we combine Eq. \ref{eqn:general-g2-einstein}, the evolution equation
\begin{equation}
\label{eqn:general-system}
    x_{k+1}(\xi) = \varphi(x_{k}(\xi)),
\end{equation}
and apply the chain rule. The 1D case, however, was computationally simpler because the tangent equations for $a$ and $q$ were regularized by the scalar $\alpha$ every time step preventing the tangent solutions from blow-ups due to the positive LE (i.e., the butterfly effect). Note that here we need to compute all first- and second-order parametric derivatives of the chart to compute $g$. Since we strive to derive a recursive relation, we regularize tangent equations in a fashion analogous to the approach in Section \ref{sec:1D}. To achieve this goal, one can recursively orthonormalize the parametric gradient through an iterative linear transformation of the parameterization and fixing $\xi=0$. In particular, we change variables from step $k$ to $k+1$ such that $$\xi_{k+1}=R_{k+1}\,\xi_{k}.$$
Note that at $\xi = 0$ we stay on the same trajectory despite the transformation. This particular choice of $\xi$ does not restrict our algorithm to concrete trajectories. Indeed, we want to ``visit" all infinitesimally small $\mu$-typical regions of the attractor after an infinite number of time steps, regardless of the choice of the initial condition. Therefore, we can always linearly re-scale the feasible space of $\xi$ such that $\xi = 0$ for our arbitrary choice of the initial condition. To simplify the notation, we skip the argument in our notation whenever $\xi=0$; for example, we use the short-hand notation $x_k(0):=x_k$, $Q(x_k(0)):=Q_k$, and so forth. Thanks to this particular transformation, the parametric gradient is automatically orthonormalized, because the chain rule implies that
\begin{equation}
\label{eqn:general-orth2}
    \nabla_{\xi_{k+1}}x_{k+1} = \nabla_{\xi_{k}}x_{k+1}\;R_{k+1}^{-1} = Q_{k+1},
\end{equation}
or, equivalently,
\begin{equation}
\label{eqn:general-orth3}
    \partial_{\xi_{k+1}^{(i)}}x_{k+1} = \partial_{\xi_k^{(j)}}x_{k+1}\;(R_{k+1}^{-1})^{(ji)} =  Q^{(:i)}_{k+1}.
\end{equation}
It means that the parametric gradient of the chart has an orthonormal basis of the column space in the updated coordinate system. Note the $R$ matrix represents the Jacobian of the step-to-step parametric transformation, i.e., $R_{k+1} = \partial\xi_{k+1}/\partial\xi_{k}$. In an analogous manner, we can derive a similar relation for the Hessian of $x_{k+1}$, represented by an $n\times m\times m$ tensor,
\begin{equation}
\label{eqn:general-orth5}
    \partial_{\xi_{k+1}^{(i)}}\partial_{\xi_{k+1}^{(j)}}x_{k+1} = \partial_{\xi_{k}^{(p)}}\partial_{\xi_{k}^{(q)}}x_{k+1}\;(R_{k+1}^{-1})^{(pi)}\;(R_{k+1}^{-1})^{(qj)}.
\end{equation}
The major benefit of the variable change is a dramatic simplification of Eq. \ref{eqn:general-g2-einstein}. Indeed, in the orthonormalized system, the $R$ matrix reduces to the identity matrix, while the norm of each column of the parametric gradient equals 1. This gives rise to the following expression for $g$,
\begin{equation}
    \label{eqn:general-g3}
    \begin{split}
    g^{(i)}_{k+1} = & - \mathrm{tr}\left(Q^T_{k+1}\; \partial_{\xi_k^{(i)}}\nabla_{\xi_{k}}x_{k+1}\right) = \\& - Q^{(:j)}_{k+1}\cdot \partial_{\xi^{(i)}_{k+1}}\partial_{\xi^{(j)}_{k+1}}x_{k+1} := - Q^{(:j)}_{k+1}\cdot a^{(i,j)}_{k+1}.
    \end{split}
\end{equation}
Consequently, only two ingredients are necessary to compute the density gradient function at $\xi = 0$. First, we need the orthogonal basis of the column space of the parametric gradient $\nabla_{\xi_k}x_{k+1}$. A recursive formula for the basis can be obtained by differentiating the system in Eq. \ref{eqn:general-system} and performing QR factorization, i.e.,
\begin{equation}
    \label{eqn:general-QT}
    \nabla_{\xi_k}x_{k+1} = D\varphi_k\,\nabla_{\xi_k}x_{k} = Q_{k+1}\,R_{k+1}.
\end{equation}
Therefore, per Eq. \ref{eqn:general-orth2}, we automatically obtain the orthonormal parametric gradient at the next time step without the need of inverting $R_{k+1}$. Since the orthonormalization is performed in a recursive manner, $\nabla_{\xi_k}x_{k} = Q_k$ by construction. To complete the algorithm, we also need a recursion for $a$. This equation can be naturally derived by differentiating Eq. \ref{eqn:general-system} twice, which gives rise to
\begin{equation}
    \label{eqn:general-AT}
\begin{split}
    \partial_{\xi_k^{(i)}}\partial_{\xi_k^{(j)}}x_{k+1} &= D^2\varphi_k(\partial_{\xi_k^{(i)}}x_k, \partial_{\xi_k^{(j)}}x_k) + D\varphi_k\,\partial_{\xi_k^{(i)}}\partial_{\xi_k^{j}}x_{k} \\ &
    =D^2\varphi_k(Q_k^{:i}, Q_k^{:j}) + D\varphi_k\,a_k^{(i,j)}.
\end{split}
\end{equation}
Note that in order to compute the SRB density gradient at step $k+1$, we need to apply the Hessian re-scaling described by Eq. \ref{eqn:general-orth5} to retrieve $a_{k+1}$. We summarize this algorithm and carefully analyse its computational properties in Section \ref{sec:general-algorithm}.

The procedure in Eq. \ref{eqn:general-QT} reduces to the recursion in Eq. \ref{eqn:1D-qk} if $m = 1$. Regardless of the choice of initial condition $Q_0$, the column vectors of $Q_k$ rigorously converge to backward Lyapunov vectors as $k\to\infty$ \cite{kuptsov-lyapunov}. The set of these column vectors is in fact an orthonormal basis of the unstable (expanding) subspace $E^u_k$ of the tangent space $TM_k$. Specific directions of backward Lyapunov vectors, however, depend on the choice of $Q_0$. Therefore, in this case, the ``convergence" should be understood that, for any orthonormal $Q_0$, the column space of $Q_k$ is guaranteed to coincide with some orthonormal basis of $E^u_k$ if $k\to\infty$. A similar procedure can be used to compute all $n$ Lyapunov vectors, including those corresponding to the negative LEs, spanning the stable (contracting) subspace $E^s_k$. In uniformly hyperbolic systems, $TM_k=E^u_k\oplus E^s_k$ at every $k$, and both the subspaces are $D\varphi$-invariant (or {\it covariant}). The covariance property implies that the product $D\varphi_k Q_k$, which we compute in Eq. \ref{eqn:general-QT}, outputs $m$ vectors that belong to the unstable subspace of the tangent space at the next time step, $TM_{k+1}$. In general, the new vectors are not orthonormal. By performing the QR factorization, however, we obtain an orthonormal basis of the unstable subspace at $k+1$. Therefore, the components of $R_{k+1}$ contain projections of the column vectors of $D\varphi_k Q_k$ onto the basis vectors of $E^u_{k+1}$.

We also observe the general recursion for the acceleration vector $a$ in Eq. \ref{eqn:general-AT} can be simplified to its one-dimensional counterpart in Eq. \ref{eqn:1D-ak} if $m=1$. Using the properties of uniform hyperbolicity, the authors of \cite{chandramoorthy-s3-new,chandramoorthy-clv} rigorously show the recursion for $a$ (a.k.a. unstable manifold curvature equation) in systems with one positive LE rigorously converges at an exponential rate. To the best of our knowledge, no rigorous results for higher-dimensional cases exist. The proof of convergence for systems with one-dimensional unstable subspaces uses the fact $a_k$ can be expressed as $C_k + D\varphi_{k-1}...D\varphi_{0}\,a_0/\prod_{i=0}^{k-1}\alpha_{k-1}^2$, where $C_k$ does not depend on $a_0$ (see Section \ref{sec:1D} for the notation explanation). By the uniform expansion property, the $a_0$-dependent term exponentially converges to zero if $k\to\infty$. In case of general systems, we find similar dependencies between $a_k^{i,j}$ and all initial conditions for the second-order tangent equation. Here, instead of re-scaling with respect to the length of the projection of $D\varphi_k\,q_k$ onto $q_{k+1}$, we are re-scaling with respect to the collection of projections included in the $R_{k+1}$ matrix. Moveover, the process of computing $a$ in the general case involves inverting $R_{k+1}$, not just a scalar, which makes it hard to directly apply the properties resulting from the uniform hyperbolicity assumption. Therefore, in this paper, we resort to an empirical study of the convergence of our algorithm (see Section \ref{sec:general-algorithm}). Note that even if the recursion converges, the specific direction of $a$ is not unique at any point on $M$, because $Q$ is also not unique as discussed above. Their product, however, that equals the SRB density gradient $g$ is unique by construction.

\subsection{General algorithm for high-dimensional systems}\label{sec:general-algorithm}
We provide a practicable algorithm based on the derivation presented in the previous section. In addition, we carefully analyse its computational cost, memory requirements, and numerically investigate its convergence. Algorithm \ref{alg1} summarizes all steps necessary to numerically compute the SRB density gradient at $N$ points along a trajectory initiated at $x_{0}\in M\subset\mathbb{R}^n$. The only optional step is included in Line 1; this step is meant to compute the dimension of the unstable subspace/manifold $m$. For many chaotic maps, this parameter is known {\it a priori} and therefore Line 1 can be skipped. If this is not the case, however, one can apply Benettin et al.'s numerical procedure \cite{benettin-les} to approximate a subset of the spectrum of Lyapunov exponents. This procedure requires solving $i\in\mathbb{Z}^{+}$ homogeneous tangent equations to identity $i$ largest LEs. The parameter $T$ represents the trajectory length and affects the accuracy of LE approximation. If the LE spectrum is evidently separated from the origin (i.e., the value of 0), then $T$ does not need to be large. Lines 3-22 of Algorithm \ref{alg1} represent the main time for-loop that computes the $g$ vector at one point on the manifold per iteration. Inside this loop, we distinguish five major stages: 1) advancing first-order tangent equation and QR factorization (Eq. \ref{eqn:general-QT}), 2) advancing second order tangent equations (Eq. \ref{eqn:general-AT}), 3) inverting the $R$ matrix and rescaling the acceleration vector $a$ (Eq. \ref{eqn:general-orth5}), 4) evaluating $g$ (Eq. \ref{eqn:general-g3}), and 5) transitioning to the next time step; updating the Jacobian and Hessian.    

\begin{algorithm}\label{alg1}
\SetAlgoLined
\SetKwInOut{Input}{Input}
\SetKwInOut{Output}{Output}
\Input{$N$, $T$, $x_0$, $n = \mathrm{size}(x_0)$}

$m = \mathrm{Benettin}(T)$ \textbf{if} $m$ unknown\;
Randomly generate $Q_0$, $a^{(i,j)}_0$ such that $\mathrm{ncol}(Q_0)=m$, $\mathrm{nrow}(Q_0)=\mathrm{size}(a_0^{(i,j)})=n$, $Q_0^TQ_0 = I$, and $i,j=1,...,m$\; 

\For(\tcp*[h]{main time loop}){$k = 0,...,N-1$}{

$S_k = D\varphi_k\,Q_k$\;
QR-factorize: $Q_{k+1}\,R_{k+1} = S_k$\;
Invert $R_{k+1}$\;

\For(\tcp*[h]{2nd-order tangent equations}){$i = 1,...,m$}{
\For{$j = 1,...,i$}{
$\tilde{a}_{k+1}^{(i,j)} = D^2\varphi_k(Q_k^{(:i)}, Q_k^{(:j)}) + D\varphi_k\,a_{k}^{(i,j)}$\;
}
}

\For(\tcp*[h]{re-scaling}){$i = 1,...,m$}{
\For{$j = 1,...,i$}{
$a_{k+1}^{(i,j)} = \tilde{a}_{k+1}^{(p,q)}\,(R^{-1})_{k+1}^{(pi)}\,(R^{-1})_{k+1}^{(qj)}$\;
}
}

\For(\tcp*[h]{evaluating $g$}){$i = 1,...,m$}{
$g_{k+1}^{(i)} = -Q_{k+1}^{(:j)}\cdot a_{k+1}^{(i,j)}$\;
}

$x_{k+1} = \varphi(x_k)$\;
Evaluate: $D\varphi_{k+1}$ and $D^2\varphi_{k+1}$\;
}
\Output{$g^{(i)}_k$, $i=1,...,m$, $k=1,...,N-1$}
\caption{SRB density gradient}
\end{algorithm}
Table \ref{tab1} summarizes the computational cost of Algorithm 1. The third column of this table includes the number of the floating point operations required in each stage as a function of the trajectory length ($N$ or $T$), system dimension $n$, and unstable manifold dimension $m$. Note the third column includes only the leading term of the flop count. The final two stages involve evaluations of nonlinear equations defined by $\varphi$ and thus their computational cost is problem-dependent. In many physics-inspired chaotic systems the cost of Lines 20-21 is relatively low. Consider the Lorenz '63 system, for example. In this case, we can think of $\varphi$ as a time discretization operator of the continuous-in-time system. For Lorenz '63, the Jacobian $D\varphi$ involves a collection of linear terms proportional to the coordinates of $x$, while the Hessian $D^2\varphi$ is constant. In many scientific/engineering applications, PDE models are discretized in space using schemes with local support (such as the finite element method), which implies the resulting Jacobians and Hessians of the fully-discretized system are sparse. Therefore, in these special cases, the cost of the most expensive stage of Algorithm \ref{alg1}, which involves second-order tangent equations, can potentially be reduced to $N\,n\,m^2$. Table \ref{tab1}, however, reflects the worst-case scenario in which no sparsity patterns occur. We also highlight the fact that in many high-dimensional chaotic systems $m\ll n$ \cite{blonigan-phdthesis}. Thus, if $n$ is large, the re-scaling stage (Lines 12-16) is rather cheaper than the second-order tangent equation stage (Lines 7-11).

We conclude that the leading term of the total flop count of Algorithm \ref{alg1} is proportional to $N\,n^3\,m^2$ in a general chaotic system. In many real-world problems, however, the final cost can be significantly reduced if one takes the advantage of a system's special structure. Our algorithm is moderately cheap in terms of the memory requirements. The most memory-consuming structure is the Hessian which, in the worst-case scenario, requires storing $n^3$ floats. As we pointed out above, however, in practical high-dimensional models, the actual ``size" of the Hessian might be dramatically smaller. Note also that, in order to advance tangent equations, we need to store $m$ $n$-dimensional basis vectors (i.e., column vectors of $Q$) and $\sim 1/2\,m^2$ acceleration vectors. The $1/2$ factor is a consequence of the assumed smoothness of the coordinate chart, which implies $a^{(i,j)} = a^{(j,i)}$ everywhere on the manifold. Notice also that the our procedure is in fact a one-step method, which means that all quantities at step $k+1$ require data only from step $k$. We do not need to store data generated at previous time steps.        

\begin{table}\label{tab1}
\caption{Computational cost of Algorithm 1.}
\begin{center}
\begin{tabular}{|c|c|c|} 
 \hline
 Stage Name & Line No. & Total Cost \\\hline\hline
 Computing $m$ (Benettin's algorithm) & 1 & $T\,n^2 m$\\\hline
 Generating initial conditions & 2 & $-$\\\hline 
 Advancing first-order tangent equations & 4 & $N\,n^2\,m$\\\hline 
 QR factorization (Householder) & 5 & $N\,n\,m^2$\\\hline 
 Inverting $R$ & 6 & $N\,m^3$\\\hline 
  Advancing second-order tangent equations & 7-11 & $N\,n^3\,m^2$\\\hline 
 Re-scaling $a$ & 12-16 & $N\,n\,m^4$\\\hline 
 Computing $g$ & 17-19 & $N\,n\,m^2$\\\hline 
 Advancing primal equation & 20 & Varies\\\hline 
 Evaluating Jacobian and Hessian & 21 & Varies\\
 \hline
\end{tabular}
\end{center}
\end{table}

The final aspect of our algorithm is its convergence. As we mentioned in the previous sections, the convergence of the Lyapunov vector equation and the second-order tangent (unstable manifold curvature) equation is rigorously guaranteed in uniformly hyperbolic systems if $m=1$. Moreover, the convergence rate is exponential in such systems. Since we lack generalization of these theoretical studies, we perform a numerical test to investigate convergence properties of Algorithm \ref{alg1}. For this purpose, we use the Baker's map introduced in Eq. \ref{eqn:straight-bakers}, as well as its 3D version $\varphi:[0,2\pi]^3\to[0,2\pi]^3$ defined as
\begin{equation}
\begin{split}
    \label{eqn:general-bakers}
    x_{k+1} = \varphi(x_k) = & \left(\begin{bmatrix}2\,x_{k}^{(1)} \\ 3\,x_{k}^{(2)} \\ x_{k}^{(3)}/6 + \,\pi\lfloor x_k^{(1)}/\pi\rfloor + \pi/3\lfloor x_k^{(2)}/(2\pi/3)\rfloor\end{bmatrix}\right. \\
     + & \left.\begin{bmatrix}
     s_1\, \sin(2x_k^{(1)})\,\sin(3/2\,x_k^{(2)}) \\
     s_2\, \sin(x_k^{(1)})\,\sin(3\,x_k^{(2)}) \\
     s_3\, \sin(6\,x_k^{(3)})\end{bmatrix}\right)\,\mathrm{mod}\,2\pi,
\end{split}
\end{equation}
which we shall refer to as the 3D Baker's map. This is an invertible chaotic map with two positive and one negative LEs, and seemingly hyperbolic behavior (see Appendix \ref{app:hyperbolicity} for more details). This map has two expanding directions, along the $x^{(1)}$ and $x^{(2)}$ axes, and one contracting direction along the third axis. Analogously to its 2D counterpart, this map models the kneading operation. The dough is extended by the factor 2 and 3 along the two orthogonal directions on the $x^{(1)}\--x^{(2)}$ plane, cut into $2\cdot 3 = 6$ squares, which are subsequently stacked in the order defined by the floor functions. These history-dependent floor functions are used to guarantee the invertibility of the nonlinear transformation by periodically distributing the third component of $x_{k+1}$ across $[0,2\pi]$. Higher-dimensional Baker's maps have been widely used in image encryption as a convenient generalization of Bernoulli shifts \cite{pichler-bakers,mao-bakers}.   

\begin{figure}
    \centering
    \includegraphics[width = 0.66\textwidth]{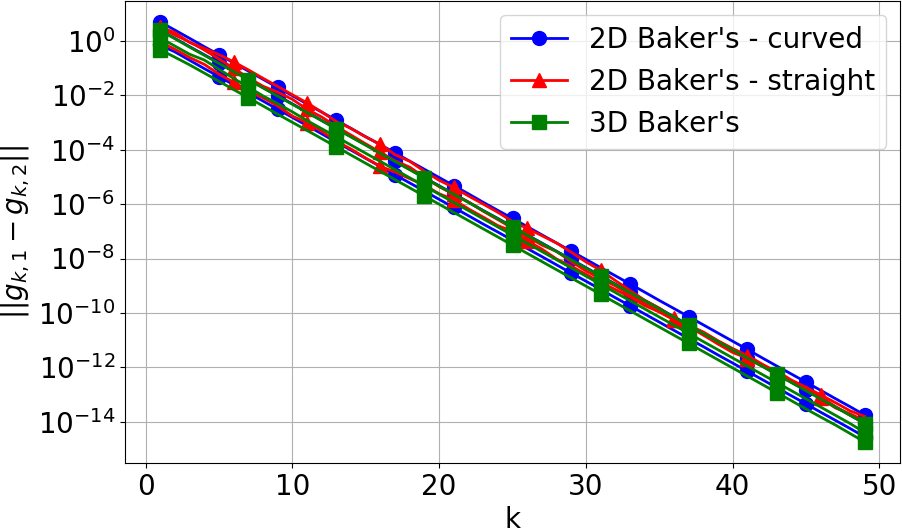}
    \caption{Relation between $\|g_{k,1}-g_{k,2}\|$ and the time step $k$ in the semilogarithmic scaling. This plot contains nine curves of three different colors. Each color corresponds to a different map: 2D Baker's map with curved unstable manifolds (blue), 2D Baker's map with straight unstable manifolds (red), and 3D Baker's map with $s_1 = 0$, $s_2 = 0.9$, $s_3 = 0.1$ (green). In case of the 2D Baker's map, the parameter values are the same as those in Figure \ref{fig:bakers_straight_pic}.}
    \label{fig:convergence}
\end{figure}

To analyse the convergence, we generate three sufficiently long trajectories started at randomly chosen initial conditions $x_0$. For each of these trajectories, we run two independent simulations with different, randomly chosen initial conditions for the tangent equations (see Line 2 of Algorithm \ref{alg1}).
Motivated by the rigorous studies, we investigate if (and how) the difference between the SRB density gradients computed along a single trajectory but using different initial conditions for tangent equations decreases in norm as we advance the iteration. In particular, we compute $\|g_{k,1}-g_{k,2}\|$, $k = 0,1,2,...$ for two random initial condition choices for tangent equations per trajectory, labelled as $1$ and $2$. The relation between this norm and time step $k$ for three different chaotic models is plotted in Figure \ref{fig:convergence}. The $g$ function is generated using Algorithm \ref{alg1} for the 2D Baker's map with $m=1$ (Eq. \ref{eqn:straight-bakers}), as well as the 3D Baker's map with $m=2$ (Eq. \ref{eqn:general-bakers}).        

We observe the norm-versus-$k$ relation is linear in the semilogarithmic scale, which clearly indicates an exponential convergence of our algorithm if applied to the Baker's map. This result implies that a relatively small number of steps ($k\approx 50$) is required to obtain the machine-precision value of the norm. Note also that the choice of trajectory ($x_0$) or model has a negligible effect on the error.

\subsection{Numerical example: Monte Carlo integration}\label{sec:general-example}

To validate Algorithm \ref{alg1}, we consider a square-integrable function $f(x)\in L^2(\mu)$ and integrate it with respect to the SRB measure $\mu$ using a Monte Carlo procedure. By the Central Limit Theorem, this integral can be approximated by taking the average of the sample distributed according to $\mu$, while the approximation error is upper-bounded by $\sqrt{\mathrm{Var}(f)/N}$, i.e.,
\begin{equation}
    \label{eqn:general-monte-carlo}
    \left|\int_{M}f(x)\,d\mu(x)-\frac{\omega(M)}{N}\sum_{k=0}^{N-1}f(x_k)\right|\leq C\sqrt{\frac{\mathrm{Var}(f)}{N}},
\end{equation}
where $C>0$ and $x_{k+1}=\varphi(x_k)\in M$. Therefore, by generating a sufficiently long trajectory and evaluating $f$ at every point along it, we gradually approach the sought-after solution. Motivated by particular applications of the SRB density gradient function (see Sections \ref{sec:introduction}-\ref{sec:preliminaries}), we consider $f(x):=\sum_{j=1}^{m}\partial_{Q^{(:j)}}v(x)$, where $v(x):M\to\mathbb{R}$ is some smooth function. In other words, we strive to integrate a sum of $m$ directional derivatives along $m$-dimensional unstable manifolds of the scalar function $v(x)$. Note integrals of this type are critical in the sensitivity computation using, for example, the general S3 method \cite{chandramoorthy-s3,chandramoorthy-s3-new}. Thanks to the partial integration (see Eq. \ref{eqn:pre-integral-lhs}-\ref{eqn:pre-integral-rhs}), we can apply Monte Carlo to two alternative versions of the same integral, since
\begin{equation}
    \label{eqn:general-integration}
    \int_M \sum_{j=1}^{m}\partial_{Q^{(:j)}}v(x)\,d\mu(x) = I = - \int_M \sum_{j=0}^m g^{(j)}(x)\,v(x)\,d\mu(x). 
\end{equation}
%Add comment about the boundary term (optional)
Using this equation, we validate Algorithm \ref{alg1} for $g$ by comparing numerical approximations of the LHS and RHS. Due to its trajectory-driven structure, Algorithm \ref{alg1} is naturally compatible with the Monte Carlo procedure. 

Two different maps shall be tested. First, we shall consider the 2D Baker's map (Eq. \ref{eqn:straight-bakers}) with $s_{4} = 0.4$ and $s_{1}=s_{2}=s_{3}=0$. As illustrated in Figure \ref{fig:bakers_straight_pic}, its unstable manifolds are curved and therefore the simplified version of the recursion for $g$ (Eq. \ref{eqn:straight-simple}) cannot be used. In this particular case, $q$ has in fact two nonzero components. Indeed, we numerically estimate that $$\max_{k\in{1,2,...,N}} \arctan\left|\frac{q_k^{(2)}}{q_k^{(1)}}\right| \approx 0.24\,\mathrm{rad} \approx 14^{\circ},$$ which is consistent with the illustration of unstable manifolds in Figure \ref{fig:bakers_straight_pic}. The second map is the 3D Baker's map (Eq. \ref{eqn:general-bakers}) with $s_{1}=0, s_{2}=0.9$, $s_{3}=0.1$. One can easily verify unstable manifolds of this map are flat surfaces aligned with the $x^{(1)}\text{--}x^{(2)}$ plane. These expanding surfaces could be curved by adding an $x^{(3)}$-dependent perturbation term to the third component of the map. 

\begin{figure}
    \centering
    \includegraphics[scale = 0.38]{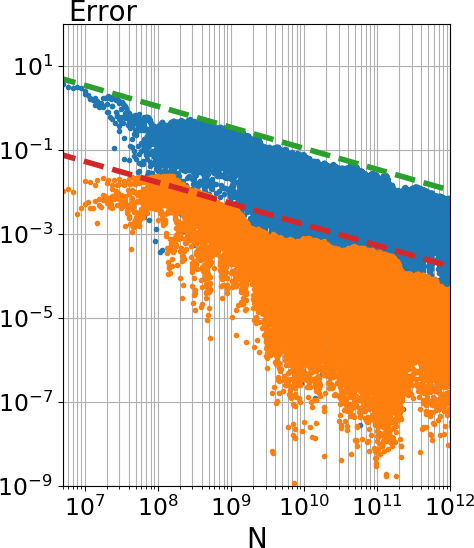}
    \includegraphics[scale = 0.38]{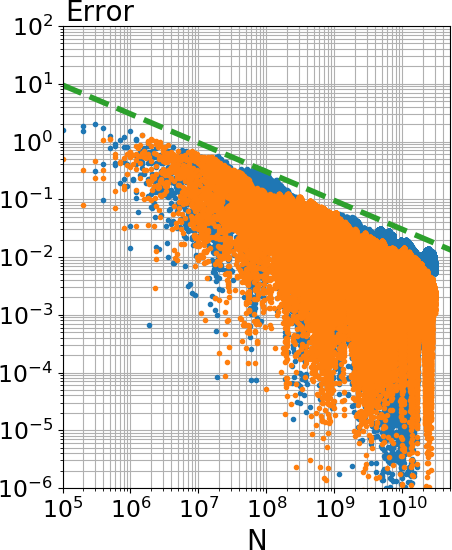}
    \caption{Error of the Monte Carlo approximation of the LHS (blue dots) and RHS (orange dots) of Eq. \ref{eqn:general-integration} versus the amount of data $N$ used. Left: 2D Baker's map. Here, we compute the relative error with respect to the reference value -1.05335809 (which equals the approximation of the RHS integral at $N=10^{13}$) for $v(x) = \sin(x^{(1)})\,\exp(x^{(2)})$. Right: 3D Baker's map. Here, we compute the absolute error with respect to the reference value of 0 for $v(x) = \sin(x^{(1)})\,\sin(3/2\,x^{(2)})\,x^{(3)}$. The dashed lines represent the slope $-1/2$ in the logarithmic scaling.}
    \label{fig:bakers_curved_convergence}
\end{figure}

Figure \ref{fig:bakers_curved_convergence} includes results of the integration test.  
Our primary conclusion is that the Monte Carlo approximations of the LHS and RHS of Eq. \ref{eqn:general-integration} approach each other as $N\to\infty$ with the rate $\mathcal{O}(1/\sqrt{N})$, which directly confirms the correctness of Algorithm \ref{alg1}. Recall we require $g$ to regularize the linear response formula, as it involves derivatives of strongly-oscillatory functions (see Section \ref{sec:preliminaries}). The examples presented in this section, however, include mildly-oscillatory functions $v(x)$ with derivatives that behave similarly (note they involve a combination of trigonometric, exponential and linear functions). Nevertheless, we observe significantly smaller errors of the RHS approximation in the 2D Baker's map case. Note the approximation error of Monte Carlo integration also depends on the variance of the integrand, which can be upperbounded by a quantity proportional to the $L^2(\mu)$-norm of the SRB density gradient $g$, denoted by $\|g\|_{L^2(\mu)}$. Indeed, $\|g\|_{L^2(\mu)}$ equals $\mathcal{O}(10^{-2})$ and $\mathcal{O}(10^{1})$ for the 2D and 3D Baker's map, respectively. This explains the significantly better performance of the Monte Carlo procedure in the former case. Therefore, if $\|g\|_{L^2(\mu)}$ exists and is sufficiently small, Monte Carlo integration might be significantly cheaper if applied to the regularized integrals of this type, regardless of the behavior of $v(x)$. If $g$ is not even Lebesgue-integrable, i.e. $g\notin L^1(\mu)$, the integrals in Eq. \ref{eqn:general-integration} do not converge, as showed in \cite{sliwiak-differentiability}.

\section{Conclusions}\label{sec:conclusions}
Ruelle's linear response formula is fundamental in the construction of numerical methods for sensitivity analysis of $n$-dimensional hyperbolic chaotic systems. Its original form, however, is impractical for direct computation due to the presence of derivatives of composite functions that grow exponentially in time. Fortunately, it is possible to easily regularize this expression through partial integration. In case of nonuniform measures describing the statistics of chaos, the by-product of the integration by parts, per the generalized fundamental theorem of calculus, involves the SRB density gradient $g$ defined as the directional derivative of conditional SRB density on $m$-dimensional unstable manifolds. Computation of $g$ is the price that must be paid for a computable version of Ruelle's formula.  

Using the measure-based coordinate parameterization, the time evolution of the measure gradient is rigorously derived by applying the measure preservation property, differentiating the coordinate charts with the chain rule on smooth manifolds. Indeed, $g$ can be computed in a recursive manner by solving a set of $\mathcal{O}(m)$ first- and $\mathcal{O}(m^2)$ second-order tangent equations, as well as step-by-step QR-factorization and inversion of $n\times m$ and $m\times m$ matrices, respectively. While the total cost of approximating $g$ at $N$ consecutive points along a trajectory is $\mathcal{O}(Nn^3m^2)$ in the worst-case scenario, the actual computational cost may scale linearly with the dimension of the system in many real-world models due to their sparse structure. Moreover, this procedure requires storing $\mathcal{O}(m^2)$ $n$-dimensional vectors only from the current time step to advance the iteration in time. Therefore, in terms of the hardware requirements, our algorithm would definitely be a reasonable choice for high-dimensional physical systems since $m\ll n$.  

The algorithm we propose is compatible with existing methods for sensitivity analysis that stem from the linear response theory, including the space-split sensitivity (S3) and FDT-based methods. Many of them approximate sensitivities through an ergodic-averaging Monte Carlo procedure and require knowledge of the directional derivative of conditional SRB measures. Moreover, $g$ can be used to assess the differentiability of statistical quantities in hyperbolic systems, which a recurring theme in theoretical studies of chaos. Thus, we believe our method provides a new major tool for both rigorous analysis and applied studies of large chaotic systems.

\section*{Acknowledgments}
The authors acknowledge the MIT SuperCloud and Lincoln Laboratory Supercomputing Center for providing HPC resources that have contributed to the research results reported within this paper.

\appendix

\section{Applying the simplified recursive formula for SRB density gradient to 1D non-injective maps}\label{app:noninjective}

Throughout this paper, we assume $\varphi$ is an invertible map. Based on this assumption, we directly use the measure preservation property to derive a recursive formula for $g$, including the simplified version for maps with straight unstable manifolds, as described in Section \ref{sec:1D-example1}. However, in the literature, one can find several one-dimensional maps such as the sawtooth/Bernoulli map \cite{sliwiak-1d}, cusp map \cite{mehta-cusp}, logistic map \cite{wormell-logistic}, onion map \cite{sliwiak-differentiability}, tent map \cite{baladi-tent}, and so forth. All of them are scientifically relevant, as they represent some simplified physics or feature interesting mathematical properties. However, most of them are non-injective, which violates the basic assumption of our derivation. In this section, however, we argue that Eq. \ref{eqn:straight-simple} can still be used to compute $g$ for such maps.

Many of the popular 1D chaotic maps (such as those listed above) are two-to-one. Thus, we assume $\varphi$ satisfies this condition; however, the argument we present can be naturally extended to other types of surjection. Let us also assume, without loss of generality, $\varphi:[0,1]\to[0,1]$ and $\varphi$ is monotonic in $[0,0.5)$ and $(0.5,1]$. Let us now define a two-dimensional analog of $\varphi$, denoted by $\varphi_{2D}:[0,1]^2\to[0,1]^2$ and satisfying
\begin{equation}
    \label{eqn:appA-1}
    x_{k+1} = \varphi_{2D}(x_k) = \begin{bmatrix}\varphi(x_k^{(1)}) \\x_k^{(2)}/2 + 0.5\lfloor2x_k^{(1)}\rfloor \end{bmatrix}.
\end{equation}
Note $\varphi_{2D}$ is invertible and resembles the 2D Baker's map (see Eq \ref{eqn:straight-bakers}). The invertibility is guaranteed by adding the floor function in $\varphi^{(2)}_{2D}$. Analogously to the 2D/3D Baker's map, here the discontinuity point is located at $x^{(1)}=0.5$, which means that the value of 0.5 is added to $x_k^{(2)}/2$ if $x_k^{(1)}>0.5$. If the monotonicity breaking point was different, then the coefficients of the floor function would need to be modified accordingly. One of the main messages of this example is to point out that any surjective 1D map can be represented as a higher-dimensional invertible map with one positive Lyapunov exponent. 

Note 1D unstable manifolds of $\varphi_{2D}$ are aligned with the first phase space coordinate, per the argument given in Section \ref{sec:1D-example1}. Thus, its SRB distribution is similar to the one of Baker's map presented in Figure \ref{fig:bakers_straight_pic}. Note also that the horizontal deformation of the trajectory of $\varphi_{2D}$ is determined solely by $\varphi$. This implies that the SRB distribution of $\varphi$ is in fact an integral of SRB distributions of $\varphi_{2D}$ restricted to single unstable manifolds over all values of $x^{(2)}$. In other words, $\varphi_{2D}$ scatters the SRB measure of $\varphi$ (which is supported on [0,1]) over an infinite set of vertically stacked intervals $[0,1]$ (which geometrically coincide with unstable manifolds of $\varphi_{2D}$). This further implies the SRB density of $\varphi$ equals the SRB distribution of $\varphi_{2D}$ integrated with respect to the vertical (second) coordinate. 

In case of a map defined by Eq. \ref{eqn:appA-1}, the simplified recursive formula for $g$ can be expressed in terms of phase space derivatives of $\varphi$ (see Section \ref{sec:1D-example1} for the derivation),
\begin{equation}
    \label{eqn:appA-2}
    g(\varphi(x)) = \frac{g(x)}{\varphi'(x)} - \frac{\varphi''(x)}{\varphi'(x)^2}.
\end{equation}
Here, the prime symbol ($'$) denotes differentiation with respect to phase space. Let $\rho(x)$ be the SRB density of $\varphi$. The $g(x)$ function that satisfies Eq. \ref{eqn:appA-2} is not the SRB density gradient of $\varphi$, defined as $g^{\varphi}(x):=\rho'(x)/\rho(x)$. According to our discussion in Section \ref{sec:1D-example1}, $g(x)$ is in fact a conditional SRB density gradient of $\varphi_{2D}$ associated with the unstable manifold parameterized by $x^{(2)}$. However, as we discussed in the previous paragraph, the SRB measure of $\varphi$ can be computed by integrating ``slices" of the SRB measure of $\varphi_{2D}$ parallel to $x^{(1)}$. This implies that, given the definition of the SRB density gradient, $g^{\varphi}(x)$ can be computed by Lebesgue-integrating the SRB density gradients obtained in the above iteration along the vertical axis. 

In practice, to construct a trajectory-based algorithm for $g^{\varphi}$, we can directly use the recursion in Eq. \ref{eqn:appA-2}. The algorithm we propose is the following. Divide the phase space $[0,1]$ in $K\in\mathbb{Z}^+$ bins of equal width. Generate a sufficiently long sequence $\{g_0, g_1, g_2,...\}$ using Eq. \ref{eqn:appA-2} starting from a random initial condition $g_0$. For each bin, take the average of the members of the sequence that correspond to one bin. Based on our discussion above, the obtained average value converges to $g^{\varphi}$. This algorithm in fact provides a piecewise constant approximation of $g^{\varphi}$.

\begin{figure}
\centering
    \includegraphics[width = 0.47\textwidth]{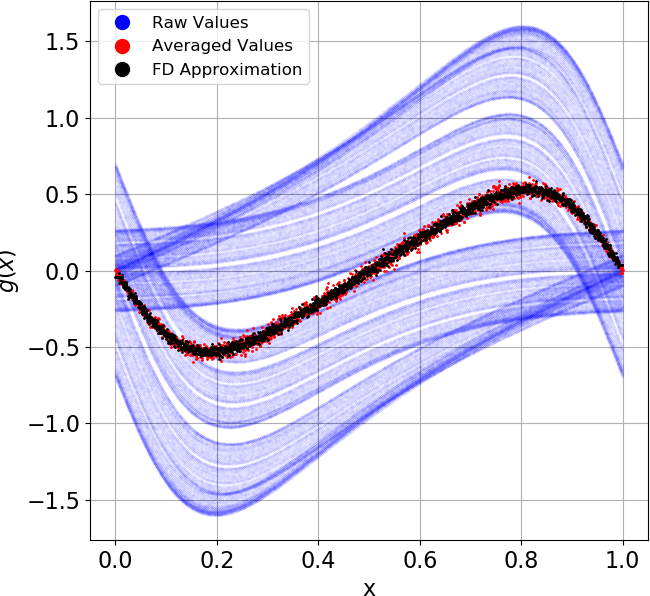}
    \includegraphics[width = 0.47\textwidth]{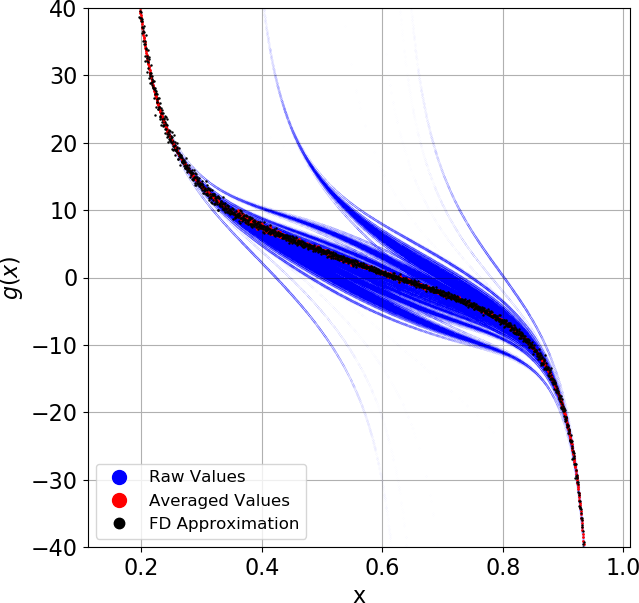}
    \caption{SRB desity gradient generated for the sawtooth map $x_{k+1} = 2x_k + s\,\sin(2\pi\,x_k)\;\mathrm{mod}\;1$ at $s = 0.1$ (left) and the onion map $x_{k+1} = 0.97\sqrt{1 - |1-2x_k|^{\gamma}}$ at $\gamma = 0.4$ (right). The averaged values (red dots) were computed by averaging the raw values (blue dots) in each of 2048 bins. The FD Approximation data points represent the central finite difference approximation of the SRB density gradient using the definition of $g$ and empirically computed SRB densities. We generated a trajectory of length $N = 10^6$ to compute the raw/averaged values of $g$.}
    \label{fig:density_gradient_1D}
\end{figure}

To verify our argument, we present a numerical experiment in which we apply the algorithm to two different 1D maps, the sawtooth map and onion map. Both of them are two-to-one and piecewise smooth. Figure \ref{fig:density_gradient_1D} shows raw values of the sequence $\{g_0, g_1, g_2, ...\}$ obtained using Eq. \ref{eqn:appA-2}, their averaged values, and finite-difference (FD) approximation of $g^{\varphi}$ using empirical SRB densities of these maps. We observe there is a good agreement between the averaged values and FD approximations in both cases. 

\begin{figure}
    \centering
    \includegraphics[width = 0.35\textwidth]{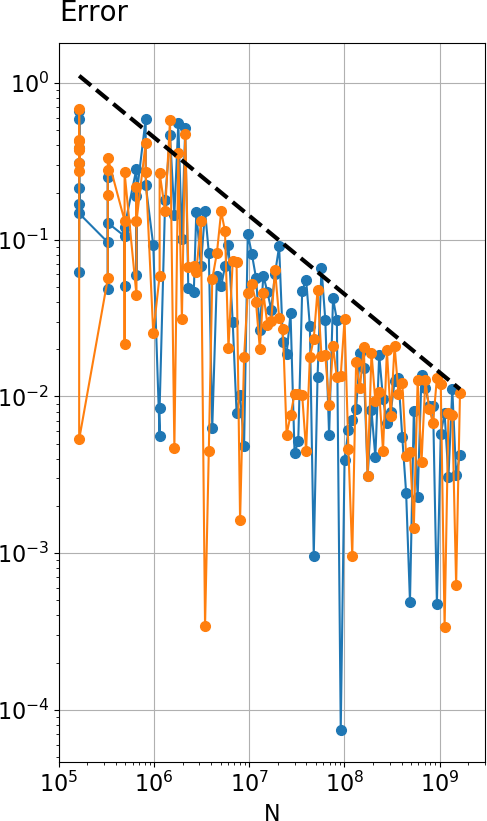}
    \caption{Relative error of the approximation of $g^{\varphi}(x)$ versus the trajectory length $N$. The error was computed for the sawtooth map at $s = 0.1$ at two phase space coordinates, $x \approx 0.4$ (blue curve) and $x \approx 0.6$ (orange curve). All error values were computed with respect to the reference value generated using $N = 3.3\cdot10^{11}$ samples. The reference dashed line represents the slope $-1/2$ in the logarithmic scaling.}
    \label{fig:convergence_sawtooth}
\end{figure}

Finally, we perform the relative error convergence test of the averaged values with respect to the trajectory length $N$. We focus on two different bins and compute the relative error with respect to a reference value generated using significantly more samples. Our results generated for the sawtooth map are shown in Figure \ref{fig:convergence_sawtooth}. As expected, the error decays and is upperbounded by $\mathcal{O}(1/\sqrt{N})$, which is a consequence of the Lebesgue-integration (or, equivalently, weighted averaging) of (conditional) SRB density gradients. This example shows that a trajectory of minimum length $N=10^9$ should be generated in order to obtain an approximation with a relative error smaller than $1\%$. 

\section{Probing the hyperbolicity of the Baker's map}\label{app:hyperbolicity}
\begin{figure}
    \centering
    \includegraphics[width = 0.6\textwidth]{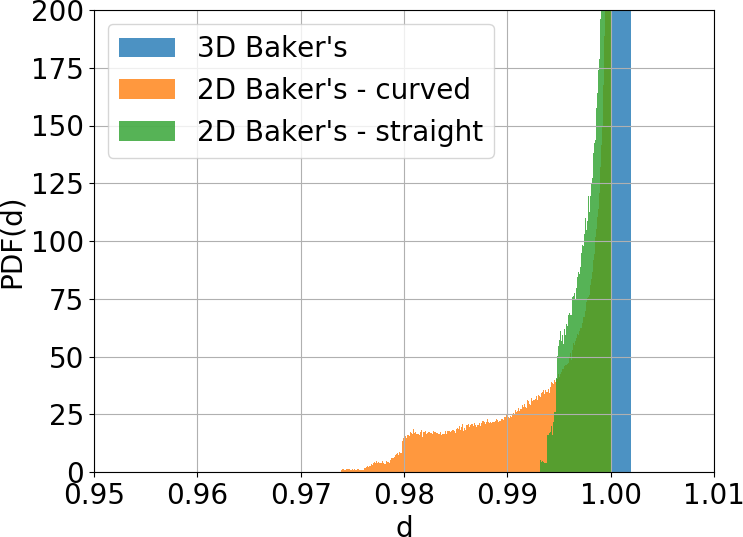}
    \caption{Hyperbolicity test performed for the 2D and 3D Bakers map. The parameter values are the same as the ones used is the numerical examples in Section \ref{sec:1D} (2D Baker's) and Section \ref{sec:general} (3D Baker's). To generate the PDF, we computed $N=10^6$ samples of $d$ along a trajectory.}
    \label{fig:hyperbolicity}
\end{figure}

Hyperbolicity guarantees the tangent space can be decomposed into two $D\varphi$-invariant subspaces, where one is asymptotically expanding (unstable), while the other one is asymptotically contracting (stable). If the expansion/contraction is uniform, then such systems are uniformly hyperbolic. Hyperbolicity is the major assumption for the dynamical systems we consider in this paper. Indeed, if the system is hyperbolic and has absolutely continuous conditional measures on unstable manifolds, then the SRB measure exists \cite{climenhaga-srb}. It is not always possible to analytically verify that a particular map is hyperbolic. Fortunately, there exist numerical procedures allowing for an efficient assessment of hyperbolicity \cite{kuptsov-lyapunov}. Most of them test the two basic criteria of hyperbolicity: 1) No zero LEs, and 2) Strict separation of the stable and unstable subspaces. Here, we apply the method proposed in \cite{kuptsov-hyperbolicity}, which computes the basis vectors of the two subspaces and approximates the smallest angle between them at different points of the manifold. If any of these angles is close to zero, then the stable and unstable subspaces are (almost) tangent, which implies the systems is likely to be non-hyperbolic. In Figure \ref{fig:hyperbolicity}, we compute the PDF of $d\in[0,1]$, which is a normalized quantity associated with the smallest principal angle between the stable and unstable subspace (our $d$ equals $k!\,d_k$; see the above reference for a rigorous definition of $d_k$). If the distribution is evidently separated from the origin ($d=0$), then it is highly likely there are no tangencies between the two subspaces. 
We observe the normalized parameter $d$ is highly unlikely to drop below the value of $0.97$. As a by-product of the applied algorithm, we computed the spectrum of Lyapunov exponents (alternatively, one can use Benettin et al.'s algorithm \cite{benettin-les}). The LEs approximately equal: $0.69\approx\log(2)$, $-0.69$ (2D Baker's with straight unstable subspaces), $0.69\approx\log(2)$, $-0.71$ (2D Baker's with curved unstable subspaces), $1.09\approx\log(3)$, $0.69\approx\log(2)$, $-1.16$ (3D Baker's). Although a small change in the parameter value does not significantly impact the LE values, it may move the PDF of $d$ closer to the origin. Based on the empirical evidence presented in this section, we conclude the 2D/3D Baker's map is clearly hyperbolic at the chosen parameter values.

\bibliographystyle{siamplain}
\bibliography{references.bib}

\end{document}